\documentclass[11pt]{article}

\usepackage{amsfonts,amsmath}
\usepackage{latexsym}

\author{K\'aroly J. B\"or\"oczky\footnote{Supported by
NKFIH grant K 132002}\footnote{Alfr\'ed R\'enyi Institute of Mathematics, 
 Realtanoda u. 13-15, H-1053 Budapest, Hungary, 
 boroczky.karoly.j@renyi.hu}, Apratim De\footnote{Department of Mathematics, Central European University, Nador u. 9, H-1051, Budapest, Hungary,
de.apratim91@gmail.com}}
\title{Stability of the
 Logarithmic Brunn-Minkowski inequality in the case of many hyperplane symmetries}

\newcommand{\proof}{\noindent{\it Proof: }}
\newcommand{\proofbox}{\mbox{ $\Box$}\\}
\newcommand{\R}{\mathbb{R}}

\newtheorem{lemma}{LEMMA}[section]
\newtheorem{theo}[lemma]{THEOREM}

\newtheorem{claim}[lemma]{CLAIM}
\newtheorem{coro}[lemma]{COROLLARY}
\newtheorem{conj}[lemma]{CONJECTURE}
\newtheorem{prop}[lemma]{PROPOSITION}

\begin{document}

\maketitle


\begin{abstract}
In the case of symmetries with respect to $n$ independent linear hyperplanes, a stability versions  of the
 Logarithmic Brunn-Minkowski Inequality and the Logarithmic Minkowski Inequality for convex bodies are established.
\end{abstract}

\noindent{\it MSC} 52A40

\section{Introduction}

The classical Brunn-Minkowski inequality forms the core of various areas in probability, additive combinatorics and convex geometry
(see Gardner \cite{Gar02}, Schneider \cite{Sch14} and Tao, Vu \cite{TaoVu}).
For recent related work in the theory of valuations, algorithmic theory and the Gaussian setting, see say
Jochemko, Sanyal \cite{JoS17,JoS18}, Kane \cite{kane}, Gardner, Zvavitch \cite{GZ},
Eskenazis, Moschidis \cite{EsM21}.
The rapidly developing new  $L_p$-Brunn-Minkowski theory (where $p=1$ is the classical case)
initiated by Lutwak \cite{Lut93,Lut93a,Lut96}, has become a main research area in modern convex geometry and geometric analysis. Following Firey \cite{Fir62} and Lutwak \cite{Lut93,Lut93a,Lut96},
major results have been obtained by
Hug, Lutwak, Yang, Zhang \cite{HLYZ2},
and more recently the papers
Kolesnikov, Milman \cite{KoM22},
Chen, Huang, Li, Liu \cite{CHLL20},
Hosle, Kolesnikov, Livshyts \cite{HKL21}, Kolesnikov, Livshyts \cite{KoL22}
 present new developments and approaches.
 We note that the $L_p$-Minkowski and $L_p$-Brunn-Minkowski inequalities are even extended to certain families of non-convex sets by
Zhang \cite{Zha99}, Ludwig, Xiao, Zhang \cite{LXZ11} and Lutwak, Yang, Zhang \cite{LYZ12}.

We call a compact convex set  $K$ in $\R^n$  a convex body if $V(K)>0$ where $V(K)$ stands for the
$n$-dimensional Lebesgue measure.
The cornerstone of the Brunn-Minkowski Theory is the Brunn-Minkowski inequality (see Schneider \cite{Sch14}).
 If $K$ and $C$ are convex bodies in $\R^n$ and $\alpha,\beta>0$, then
the Brunn-Minkowski inequality
 says that
\begin{equation}
\label{BrunnMinkowski}
V(\alpha K+\beta C)^{\frac1n}\geq \alpha V(K)^{\frac1n}+\beta V(C)^{\frac1n}
\end{equation}
where equality holds if and only if $C=\gamma K+z$ for $\gamma>0$ and $z\in\R^n$. Because of the homogeneity of Lebesgue measure, \eqref{BrunnMinkowski} is equivalent
to saying that if $\lambda\in(0,1)$, then
\begin{equation}
\label{BrunnMinkowskiprod}
V((1-\lambda) K+\lambda\, C)\geq  V(K)^{1-\lambda} V(C)^{\lambda}
\end{equation}
where equality holds if and only if $K$ and $C$ are translates. We also note another consequence of
the Brunn-Minkowski inequality \eqref{BrunnMinkowski}; namely, the Minkowski inequality says that
\begin{equation}
\label{Minkowski}
\int_{S^{n-1}} h_C\,dS_K\geq \int_{S^{n-1}} h_K\,dS_K\mbox{ \ \ provided $V(C)=V(K)$}
\end{equation}
where $S_K$ is the surface area measure on $S^{n-1}$ and $h_K(u)=\max_{x\in K}\langle u,x\rangle$ is the support function of $K$ (see Schneider \cite{Sch14}).

The first stability forms of the Brunn-Minkowski inequality for convex bodies were due to Minkowski himself
(see Groemer \cite{Gro93}).
If the distance of  $K$ and $C$ is measured in terms of the so-called Hausdorff distance, then
Diskant \cite{Dis73} and Groemer \cite{Gro88} provided close to optimal stability versions
(see Groemer \cite{Gro93}). However, the natural distance is in terms of the volume of the symmetric difference, and the optimal result is due to Figalli, Maggi, Pratelli \cite{FMP09,FMP10} (see Figalli, van Huntum, Tiba \cite{FHT} for the case of general bounded measurable sets).
To define the ``homothetic distance'' $A(K,C)$
of convex bodies $K$ and $C$, let $\alpha=V(K)^{\frac{-1}n}$ and
$\beta=V(C)^{\frac{-1}n}$, and let
$$
A(K,C)=\min\left\{V\Big(\alpha K\Delta (x+\beta C)\Big):\,x\in\R^n\right\}
$$
where $K\Delta Q$ stands for the symmetric difference of $K$ and $Q$.
In addition, let $\sigma(K,C)=\max\left\{\frac{V(C)}{V(K)},\frac{V(K)}{V(C)}\right\}$.
Now Figalli, Maggi, Pratelli \cite{FMP10} proved that setting $\gamma^*=(\frac{(2-2^{\frac{n-1}{n}})^{\frac32}}{122n^7})^2$, we have
$$
V(K+C)^{\frac1n}\geq (V(K)^{\frac1n}+V(C)^{\frac1n})
\left[1+\frac{\gamma^*}{\sigma(K,C)^{\frac1n}}\cdot A(K,C)^2\right].
$$

Here the exponent $2$ of $A(K,C)^2$ is optimal ({\it cf.} Figalli, Maggi, Pratelli \cite{FMP10}).
We note that prior to \cite{FMP10}, the only known error term in the Brunn-Minkowski inequality
was of order $A(K,C)^\eta$ with $\eta\geq n$,
due to
Diskant \cite{Dis73} and  Groemer \cite{Gro88} in their work on providing stability result in terms of the Hausdorff distance
(see also Groemer \cite{Gro93}), and also to a more direct approach by
Esposito, Fusco, Trombetti \cite{EFT05}; therefore, the exponent depended significantly on $n$.

Figalli, Maggi, Pratelli \cite{FMP10} proved a factor of the form $\gamma^*(n)=cn^{-14}$ for some absolute constant $c>0$, which was improved to $c n^{-7}$ by Segal \cite{Seg12}, and subsequently to $c n^{-5.5}$ by Kolesnikov, Milman \cite{KoM22}, Theorem~12.12. 
The current best known bound for $\gamma^*(n)$ is $cn^{-5}(\log n)^{-1}$, which follows by combining the general estimate of Kolesnikov-Milman \cite{KoM22}, Theorem 12.2, with the logarithmic bound of Klartag \cite{Kla23} on the Cheeger constant of a convex body in isotropic position improving on Yuansi Chen’s work \cite{Che21} on the Kannan-Lov\'asz-Simonovits conjecture.
Harutyunyan \cite{Har18} conjectured that $\gamma^*(n)=c n^{-2}$ is the optimal order of the constant, and 
showed that it can't be of smaller order. Actually, Segal \cite{Seg12} observed that
 Dar's conjecture in \cite{Dar99} would imply that we may choose $\gamma^*(n)=c n^{-2}$
for some absolute constant $c>0$.

We note that recently, various breakthrough stability results about geometric functional inequalities have been obtained.
Fusco, Maggi, Pratelli \cite{FMP08} proved an optimal stability version of the isoperimetric inequality (whose result was extended to the Brunn-Minkowski inequality by Figalli, Maggi, Pratelli \cite{FMP09,FMP10}, see also Eldan, Klartag \cite{ElK14}). Stonger versions of the  Borell-Brascamp-Lieb inequality are provided by Ghilli, Salani \cite{GhS17} and Rossi, Salani \cite{RoS17}, and
 of the Sobolev inequality by Figalli, Zhang \cite{FiZ} (extending Bianchi, Egnell \cite{BiE91} and Figalli, Neumayer \cite{FiN19}),
Nguyen \cite{Ngu16} and Wang \cite{Wan16}, and of some related inequalities by
Caglar, Werner \cite{CaW17}. Related inequalities are verified by Colesanti \cite{Col08}, 	Colesanti, Livshyts, Marsiglietti \cite{CLM17}, P. Nayar, T. Tkocz \cite{NaT13,NaT20},
 Xi, Leng \cite{XiL16}.

In this paper, we focus on   replacing Minkowski addition with the $L_0$ sum. First, for $\lambda\in(0,1)$, the $L_0$ or logarithmic sum of two origin symmetric convex bodies $K$ and $C$ in $\R^n$ is defined by
$$
(1-\lambda)\cdot K+_0 \lambda\cdot C=\left\{x\in\R^n:\,\langle x,u\rangle\leq
h_K(u)^{1-\lambda}h_C(u)^\lambda\;\forall u\in S^{n-1}\right\}.
$$
It is linearly equivariant, as
 $A((1-\lambda)\cdot K+_0 \lambda\cdot C)=(1-\lambda)\cdot A\,K+_0 \lambda\cdot A\,C$ for $A\in{\rm GL}(n)$.
We say that the compact convex sets $K,C\subset \R^n$ are dilates if $K=\lambda C$ for $\lambda>0$. In addition, a compact convex sets $K\subset \R^n$ is centered if the centroid of $K$ is the origin; namely, if  
$\int_K x\,d\mathcal{L}^d(x)=o$ where ${\rm dim}\,K=d$ and $\mathcal{L}^d$ is the $d$-dimensional Lebesgue measure.

The following strengthening of the Brunn-Minkowski inequality for centered convex bodies
is a long-standing and highly investigated conjecture. 

\begin{conj}[Logarithmic Brunn-Minkowski conjecture]
\label{logBMconj}
If $\lambda\in(0,1)$ and $K$ and $C$ are centered convex bodies in $\R^n$, then
\begin{equation}
\label{logBMconjeq}
V((1-\lambda)\cdot K +_0 \lambda\cdot C)\geq V(K)^{1-\lambda} V(C)^\lambda,
\end{equation}
with equality if and only if $K=K_1+\ldots + K_m$ and $C=C_1+\ldots + C_m$  for centered compact convex sets
	$K_1,\ldots, K_m,C_1,\ldots,C_m$ of dimension at least one where $\sum_{i=1}^m{\rm dim}\,K_i=n$
	and $K_i$ and $C_i$ are dilates, $i=1,\ldots,m$.
\end{conj}

We note that the choice of the right translates of $K$ and $C$ is important in Conjecture~\ref{logBMconj} according to the examples by Nayar, Tkocz \cite{NaT13}. On the other hand, the following is an equivalent form of the  origin symmetric case of the Logarithmic Brunn-Minkowski conjecture for $o$-symmetric convex bodies.

The cone volume measure or $L_0$-surface area measure $V_K$ on  $S^{n-1}$, 
whose study was initiated independently by Firey \cite{Fir74}
and Gromov and Milman \cite{GromovMilman},  has become an indispensable tool in the last decades 
(see say Barthe, Gu\'{e}don, Mendelson, Naor \cite{BG}, Naor \cite{Nar07},  Paouris, Werner \cite{PW}).
If a convex body $K$ contains the origin, then its cone volume measure is
$dV_K=\frac1n\,h_K\,dS_K$ where $h_K$ is the support function of $K$ and the total measure is the volume of $K$.

Following partial and related results by Andrews \cite{And99}, Chou, Wang \cite{CW},
He, Leng, Li \cite{HLL06},
Henk, Sch\"urman, Wills \cite{HSW06}, Stancu \cite{Stancu},
Xiong \cite{Xio10}
the paper B\"or\"oczky, Lutwak, Yang, Zhang \cite{BLYZ13} characterized even cone volume measures
by the so called subspace concentration condition. As it turns out, subspace concentration condition also holds for the cone-volume measure $V_K$ if
the centroid of a general convex body $K$ is the origin (see Henk, Linke \cite{HeL14} and B\"or\"oczky, Henk \cite{BoH16,BoH17}).

\begin{conj}[Logarithmic Minkowski conjecture]
\label{logMconj}
If $K$ and $C$ are centered convex bodies in $\R^n$, then
\begin{equation}
\label{logMconjeq}
\int_{S^{n-1}}\log \frac{h_C}{h_K}\,dV_K\geq \frac{V(K)}n\log\frac{V(C)}{V(K)}
\end{equation}
with the same equality conditions as in Conjecture~\ref{logBMconj}.
\end{conj}

In $\R^2$, Conjecture~\ref{logBMconj} and Conjecture~\ref{logMconj} are verified by B\"or\"oczky, Lutwak, Yang, Zhang \cite{BLYZ12} for $o$-symmetric convex bodies, but it is still open in general. On the other hand, Xi, Leng \cite{XiL16} proved that any two dimensional convex bodies $K$ and $C$
in $\R^2$ can be translated in a way such that \eqref{logBMconjeq} and \eqref{logMconjeq} hold for the translates.
In higher dimensions, Conjecture~\ref{logBMconj} and Conjecture~\ref{logMconj} are proved
for complex bodies ({\it cf.} Rotem \cite{Rotem}) and for convex bodies invariant that are invariant with respect to $n$ independent linear reflections  ({\it cf.} Theorem~\ref{logBMsym}).
 The latter type of bodies include unconditional convex bodies, which was handled earlier by Saroglou \cite{Sar15}.

In addition, Conjecture~\ref{logMconj} is verified if $C$ is origin symmetric and $K$ is a zonoid by van Handel \cite{vHa} (with equality case only clarified when $K$ has a $C^2_+$ boundary), or if $C$ is a convex body whose centroid is the origin and $K$ is a centered ellipsoid by Guan, Ni \cite{GuN17}. 
For origin symmetric $K$ and $C$, Conjecture~\ref{logMconj} is proved when $K$ is close to an ellipsoid (with equality case only clarified when $K$ has a $C^2_+$ boundary)  by a combination of the local estimates by Kolesnikov, Milman \cite{KoM22}, and the use of the continuity method in PDE by Chen, Huang, Li, Liu \cite{CHLL20}. Here closeness to an ellipsoid means that there exist some $c_n>0$ depending only on $n$ and an origin symmetric ellipsoid $E$ such that $E\subset K\subset (1+c_n)E$. Another even more recent proof of this result is due to Putterman \cite{Put21}. We note that an analogous result holds for linear images of Hausdorff neighbourhoods of $l_q$ balls for $q>2$ if the dimension $n$ is high enough according to \cite{KoM22} and the method of \cite{CHLL20}.
Actually, E. Milman \cite{Mila,Milb} (see also Ivaki,  Milman \cite{IvM23}) provides rather generous explicit curvature pinching bounds for $\partial K$ in order to Conjecture~\ref{logMconj} to hold, and proves that for any origin symmetric convex body $M$ there exists an
origin symmetric convex body $K$ with $C^\infty_+$ boundary and $M\subset K\subset 8M$ such that Conjecture~\ref{logMconj} holds for any origin symmetric convex body $C$.
Additional local versions of Conjecture~\ref{logMconj} are due to Colesanti, Livshyts, Marsiglietti \cite{CLM17}, Kolesnikov, Livshyts \cite{KoL22} and Hosle, Kolesnikov, Livshyts \cite{HKL21}.

We say that $A\in{\rm GL}(n)$ is a linear reflection associated to a linear $(n-1)$-space $H\subset \R^n$
if $A$ fixes the points of $H$ and $\det A=-1$. In this case, there exists $u\in\R^n\backslash H$ such that $Au=-u$ where
the invariant subspace $\R u$ is uniquely determined  (see Davis \cite{Dav08}, Humphreys \cite{Hum90}, Vinberg \cite{Vin71}).
It follows that a linear reflection $A$ is a classical "orthogonal" reflection if and only if $A\in O(n)$.
We say that the linear reflections $A_1,\ldots,A_n$ associated to linear $(n-1)$-spaces $H_1,\ldots,H_n\subset \R^n$ are {\it independent} if $H_1\cap \ldots\cap H_n=\{o\}$. In this case, the $n$ linear $(n-1)$-spaces $H_1,\ldots,H_n\subset \R^n$ are also called {\it independent}. We also observe that if a compact convex set $K\subset\R^n$ is invariant under $n$ independent linear reflections, then $K$ is centered.

Following the result on unconditional convex bodies by Saroglou \cite{Sar15}, B\"or\"oczky, Kalantzopoulos \cite{BoK} verified the
logarithmic Brunn-Minkowski and Minkowski conjectures
 under hyperplane symmetry assumption. 

\begin{theo}[B\"or\"oczky, Kalantzopoulos]
\label{logBMsym}
If convex bodies $K$ and $C$  in $\R^n$ are invariant under linear reflections $A_1,\ldots,A_n$ through $n$ hyperplanes $H_1,\ldots,H_n$ with $H_1\cap \ldots\cap H_n=\{o\}$, then
\begin{eqnarray}
\label{logBMsymeq}
V((1-\lambda)\cdot K +_0 \lambda\cdot C)&\geq & V(K)^{1-\lambda} V(C)^\lambda\\[1ex]
\label{logMsymeq}
\int_{S^{n-1}}\log \frac{h_C}{h_K}\,dV_K&\geq & \frac{V(K)}n\log\frac{V(C)}{V(K)},
\end{eqnarray}
with equality in either inequality
 if and only if $K=K_1+\ldots+ K_m$ and $C=C_1+\ldots + C_m$ for compact convex sets
	$K_1,\ldots, K_m,C_1,\ldots,C_m$ of dimension at least one and invariant under  $A_1,\ldots,A_n$
where $K_i$ and $C_i$ are dilates, $i=1,\ldots,m$, and $\sum_{i=1}^m{\rm dim}\,K_i=n$.
\end{theo}

Geometric inequalities under $n$ independent hyperplane symmetries were first considered by
Barthe, Fradelizi \cite{BaF13} and Barthe, Cordero-Erausquin \cite{BaC13}. These papers verified
the classical Mahler conjecture and Slicing conjecture, respectively, for these type of bodies.
The main result of our paper is a stability version of Theorem~\ref{logBMsym}.

\begin{theo}
\label{logBMstab}
If $\lambda\in [\tau,1-\tau]$ for $\tau\in(0,\frac12]$, convex bodies $K$ and $C$  in $\R^n$ are invariant under linear reflections $A_1,\ldots,A_n$ through $n$ hyperplanes $H_1,\ldots,H_n$ with
$H_1\cap \ldots\cap H_n=\{o\}$, and
$$
V((1-\lambda)\cdot K +_0 \lambda\cdot C)\leq (1+\varepsilon) V(K)^{1-\lambda} V(C)^\lambda
$$
for $\varepsilon>0$,  then for some $m\geq 1$,
there exist  compact convex sets
$K_1,C_1,\ldots,K_m,C_m$  of dimension at least one and invariant under  $A_1,\ldots,A_n$
where $K_i$ and $C_i$ are dilates, $i=1,\ldots,m$, and $\sum_{i=1}^m{\rm dim}\,K_i=n$ such that
\begin{eqnarray*}
K_1+\ldots + K_m\subset &K&\subset
\left(1+c^n\left(\frac{\varepsilon}{\tau}\right)^{\frac1{95n}}\right)(K_1+\ldots + K_m)\\
C_1+\ldots +C_m\subset &C&\subset
\left(1+c^n\left(\frac{\varepsilon}{\tau}\right)^{\frac1{95n}}\right)(C_1+\ldots +C_m)
\end{eqnarray*}
where $c>1$ is an absolute constant.
\end{theo}

Let us present an example showing that the bound of
Theorem~\ref{logBMstab} is not far from being optimal in the sense that
the exponent $1/(95n)$ should be at least $1/n$.
If for small $\varepsilon>0$, $K$ is obtained from the box $K_0=[\frac{-1}{2^{n-1}},\frac1{2^{n-1}}]\times [-2,2]^{n-1}$ by cutting off corners
of size of order $\varepsilon^{\frac1n}$, and
$C$ is obtained from the box $C_0=[-2^{n-1},2^{n-1}]\times [\frac{-1}2,\frac12]^{n-1}$ by cutting off corners
of suitable size of order $\varepsilon^{\frac1n}$, then
$\frac12\cdot K +_0 \frac12\cdot C=[-1,1]^n$, and
$$
V\left(\frac12\cdot K +_0 \frac12\cdot C\right)\leq (1+\varepsilon) V(K)^{\frac12} V(C)^{\frac12},
$$
but if $\eta K_0\subset K$ for $\eta>0$, then $\eta\leq 1-\gamma\,\varepsilon^{\frac1n}$
where $\gamma>0$ depends on $n$.

We deduce from Theorem~\ref{logBMstab} a stability version of
the logarithmic-Minkowski inequality
\eqref{logMsymeq} for convex bodies with many hyperplane symmetries.

\begin{theo}
\label{logMstab}
If  the convex bodies $K$ and $C$  in $\R^n$ are invariant under linear reflections $A_1,\ldots,A_n$ through $n$ hyperplanes $H_1,\ldots,H_n$ with
$H_1\cap \ldots\cap H_n=\{o\}$, and
$$
\int_{S^{n-1}}\log \frac{h_C}{h_K}\,\frac{dV_K}{V(K)}\leq \frac{1}n\cdot \log\frac{V(C)}{V(K)}+\varepsilon
$$
for $\varepsilon>0$,  then  for some $m\geq 1$,
there exist  compact convex sets
$K_1,C_1,\ldots,K_m,C_m$  of dimension at least one and invariant under  $A_1,\ldots,A_n$
where $K_i$ and $C_i$ are dilates, $i=1,\ldots,m$, and $\sum_{i=1}^m{\rm dim}\,K_i=n$ such that
\begin{eqnarray*}
K_1+\ldots + K_m\subset &K&\subset
\left(1+c^n\varepsilon^{\frac1{95n}}\right)(K_1+\ldots +K_m)\\
C_1+\ldots +C_m\subset &C&\subset
\left(1+c^n\varepsilon^{\frac1{95n}}\right)(C_1+\ldots +C_m)
\end{eqnarray*}
where $c>1$ is an absolute constant.
\end{theo}

If $K$ is a ball centered at the origin (and hence $m=1$), then Ivaki \cite{Iva17}, Theorem~2.1 proves an improved version of Theorem~\ref{logMstab} 
where $C$ does not need to satisfy any symmetry assumption (only translated in a suitable way) and the
error term is of order $\varepsilon^{\frac1{n+1}}$ instead of $\varepsilon^{\frac1{95n}}$.

We note that the Logarithmic Minkowski Conjecture~\ref{logMconj} is intimately related to
the Monge-Amp\`ere type logarithmic Minkowski Problem in the even case (see B\"or\"oczky, Lutwak, Yang, Zhang \cite{BLYZ12},
Kolesnikov, Milman \cite{KoM22}, Chen, Huang, Li, Liu \cite{CHLL20}, 
K.J. B\"or\"oczky \cite{Bor23}).
Recently, breakthrough
results have been obtained by
Chen, Li, Zhu \cite{CLZ19}, Chen, Huang, Li \cite{CHLL20}, Kolesnikov \cite{Kol20},
Nayar, Tkocz \cite{NaT20},
Kolesnikov, Milman \cite{KoM22},
 Putterman \cite{Put21}. Actually, Theorem~\ref{logMstab} implies the stability of the solution of the Monge-Amp\'ere equation Logarithmic-Minkowski Problem on $S^{n-1}$ for unconditional data according to B\"or\"oczky, De \cite{BoD21b}.

To prove Theorem~\ref{logBMstab}, first we verify it in the unconditional case, see Section~\ref{secuncond} presenting these partial results. More precisely, first we consider the coordinatewise product of unconditional convex bodies based on the recent stability version of the Prekopa-Leindler inequality (see Section~\ref{seccoordinatewise}), and then
handle the unconditional case Theorem~\ref{logBMuncondstab} of Theorem~\ref{logBMstab}  in
Sections~\ref{seclinimage} and \ref{seclogBMuncondstab}. Next we review some fundamental properties of Weyl chambers and Coxeter groups in general in Section~\ref{seclinear-reflections} and
Section~\ref{secCoxeter-groups}, and prove Theorem~\ref{logBMstab} in
Section~\ref{seclogBMstab}. Finally, Theorem~\ref{logMstab} is verified in
Section~\ref{seclogMstab}.

\section{The case of unconditional convex bodies}
\label{secuncond}

The way to prove Theorem~\ref{logBMstab} is to first verify the case of unconditional convex bodies; namely, when
 $A_1,\ldots,A_n$ are orthogonal reflections and $H_1,\ldots,H_n$ are coordinate hyperplanes.
For unconditional convex bodies, the coordinatewise product is a classical tool; namely, if $\lambda\in(0,1)$ and $K$ and $C$ are unconditional convex bodies  in $\R^n$, then
\begin{eqnarray*}
 K^{1-\lambda}\cdot C^\lambda&=&
\{(\pm|x_1|^{1-\lambda}|y_1|^\lambda,\ldots,\pm|x_n|^{1-\lambda}|y_n|^\lambda)\in\R^n: \\
&&(x_1,\ldots,x_n)\in K\mbox{ \ and \ }(y_1,\ldots,y_n)\in C\}.
\end{eqnarray*}
It is known (see say Saroglou \cite{Sar15}) that $K^{1-\lambda}\cdot C^\lambda$ is a convex unconditional body,
and it follows from the H\"older inequality (see also Saroglou \cite{Sar15}) that
$$
 K^{1-\lambda}\cdot C^\lambda\subset (1-\lambda)\cdot K +_0 \lambda\cdot C.
$$
In addition, \cite{Sar15} verifies that if
$\lambda\in (0,1)$, $\Phi$ is a positive definite diagonal matrix and $K$ is an unconditional convex body in $\R^n$, then
\begin{equation}
\label{KTKcoordinatewise}
K^{1-\lambda}\cdot (\Phi K)^\lambda=\Phi^\lambda K
\end{equation}
where $\Phi^\eta={\rm diag}(t_1^\eta,\ldots,t_n^\eta)$ if  $\eta\in\R$ and $\Phi={\rm diag}(t_1,\ldots,t_n)$
for $t_1,\ldots,t_n>0$.

The Logarithmic Brunn-Minkowski Conjecture~\ref{logBMconj} was verified for unconditional convex bodies
 by several authors, as
Bollobas, Leader \cite{BoL95}, 
Uhrin \cite{Uhr94} and
Cordero-Erausquin, Fradelizi, Maurey \cite{CEFM04}
verified the inequality $V((1-\lambda)\cdot K +_0 \lambda\cdot C)\geq
 V(K)^{1-\lambda} V(C)^\lambda$ in \eqref{logBMeq} about the coordinatewise product,
even before the log-Brunn-Minkowski conjecture
was stated, and the
containment relation between the coordinatewise product and the $L_0$-sum
and the description of the
equality case are due to Saroglou \cite{Sar15}. For $X,Y\subset \R^n$, we write $X\oplus Y$ to denote
$X+Y$ if ${\rm lin}X$ and ${\rm lin}Y$ are orthogonal.

\begin{theo}[Saroglou]
\label{BLCEFMS}
If $K$ and $C$ are unconditional convex bodies in $\R^n$  and $\lambda\in(0,1)$, then
\begin{equation}
	\label{logBMeq}
	V((1-\lambda)\cdot K +_0 \lambda\cdot C)\geq
V(K^{1-\lambda}\cdot C^\lambda)\geq V(K)^{1-\lambda} V(C)^\lambda.
	\end{equation}
\begin{description}	
\item[(i)] $V(K^{1-\lambda}\cdot C^\lambda)=V(K)^{1-\lambda} V(C)^\lambda$
  if and only if $C=\Phi K$ for a positive definite diagonal matrix $\Phi$.
\item[(ii)] $V((1-\lambda)\cdot K +_0 \lambda\cdot C)=V(K)^{1-\lambda} V(C)^\lambda$
  if and only if $K=K_1\oplus\ldots \oplus K_m$ and $L=L_1\oplus\ldots \oplus L_m$ for unconditional compact convex sets
	$K_1,\ldots, K_m,L_1,\ldots,L_m$ of dimension at least one where $K_i$ and $L_i$ are dilates, $i=1,\ldots,m$.
\end{description}
\end{theo}

We note that the second inequality in \eqref{logBMeq} (about the coordinatewise product) is a consequence
of the Prekopa-Leindler inequality (see  Section~\ref{seccoordinatewise}). In turn, the stability version
Theorem~\ref{PLhstablambda} of the Prekopa-Leindler inequality yields the following:

\begin{theo}
\label{coordinatewisestab}
If $\lambda\in [\tau,1-\tau]$ for $\tau\in(0,\frac12]$, and the unconditional convex bodies $K$ and $C$ in $\R^n$
satisfy
$$
V(K^{1-\lambda}\cdot C^\lambda)\leq (1+\varepsilon) V(K)^{1-\lambda} V(C)^\lambda
$$
for $\varepsilon>0$,  then
there exists  positive definite diagonal matrix $\Phi$ such that
$$
V(K\Delta (\Phi C))<c^nn^{n}\left(\frac{\varepsilon}{\tau}\right)^{\frac1{19}} V(K)
\mbox{ \ and \ }
V((\Phi^{-1} K)\Delta C)<c^nn^{n}\left(\frac{\varepsilon}{\tau}\right)^{\frac1{19}} V(C)
$$
where $c>1$ is an absolute constant.
\end{theo}

In the case of the logarithmic-Brunn-Minkowski inequality for unconditional convex bodies, we have a different type
stability  estimate:

\begin{theo}
\label{logBMuncondstab}
If $\lambda\in [\tau,1-\tau]$ for $\tau\in(0,\frac12]$, and the unconditional convex bodies $K$ and $C$ in $\R^n$
satisfy
$$
V((1-\lambda)\cdot K +_0 \lambda\cdot C)\leq (1+\varepsilon) V(K)^{1-\lambda} V(C)^\lambda
$$
for $\varepsilon>0$,  then for some $m\geq 1$,
there exist $\theta_1,\ldots,\theta_m>0$ and unconditional compact convex sets
$K_1,\ldots,K_m$ such that ${\rm lin}\,K_i$, $i=1,\ldots,m$, are complementary coordinate subspaces, and
\begin{eqnarray*}
K_1\oplus\ldots \oplus K_m\subset &K&\subset
\left(1+c^n\left(\frac{\varepsilon}{\tau}\right)^{\frac1{95n}}\right)(K_1\oplus\ldots \oplus K_m)\\
\theta_1K_1\oplus\ldots \oplus \theta_mK_m\subset &C&\subset
\left(1+c^n\left(\frac{\varepsilon}{\tau}\right)^{\frac1{95n}}\right)(\theta_1K_1\oplus\ldots \oplus \theta_mK_m)
\end{eqnarray*}
where $c>1$ is an absolute constant.
\end{theo}

\section{Coordinatewise product }
\label{seccoordinatewise}

The main tool is the  Pr\'ekopa-Leindler inequality; that is, a functional form of the Brunn-Minkowski inequality, due to
Pr\'ekopa \cite{Pre71} and Leindler \cite{Lei72} in dimension one, and to
 Pr\'ekopa \cite{Pre73}, C. Borell \cite{Bor75} and Brascamp, Lieb \cite{BrL76} in higher dimensions
 (see Artstein-Avidan, Florentin, Segal \cite{AFS20} for a recent variant).
Various applications are provided and surveyed in Ball \cite{Bal},
Barthe \cite{Bar} and Gardner \cite{Gar02}.
The following multiplicative version from \cite{Bal} is the most convenient for geometric applications.

\begin{theo}[Pr\'ekopa-Leindler]
\label{PLn}
If $\lambda\in(0,1)$ and non-negative $h,f,g\in L_1(\R^n)$
satisfy $h((1-\lambda)x+\lambda y)\geq f(x)^{1-\lambda}g(y)^\lambda$ for
$x,y\in\mathbb{R}^n$, then
\begin{equation}
\label{PLnineq}
\int_{\R^n}  h\geq \left(\int_{\R^n}f\right)^{1-\lambda} \cdot \left(\int_{\R^n}g\right)^\lambda.
\end{equation}
\end{theo}

The case of equality in Theorem~\ref{PLn} has been characterized by
 Dubuc \cite{Dub77}, and the functions $f$, $g$ and $h$ should be essentially log-concave in the case of equality.
 Here a non-negative function $\varphi$ on $\R^n$ is log-concave if
 $\varphi((1-\lambda)x+\lambda\,y)\geq \varphi(x)^{1-\lambda}\varphi(y)^{\lambda}$
 for all $x,y\in\R^n$ and $\lambda\in(0,1)$. 
In B\"or\"oczky, De \cite{BoD21a},  the following stability version of
the Prekopa-Leindler inequality for log-concave functions is verified.

 \begin{theo}
\label{PLhstablambda}
If  $\lambda\in(0,1)$ and
$f,g$ are log-concave functions on $\R^n$ satisfying $0<\int_{\R^n} f=\int_{\R^n} g<\infty$ and
$$
\int_{\R^n}\sup_{z=(1-\lambda)x+\lambda y} f(x)^{1-\lambda}g(y)^\lambda\,dz
\leq (1+\varepsilon)\int_{\R^n} f
$$
for $\varepsilon>0$, then there exists $w\in\R^n$ such that
\begin{equation}
\label{PLhstablambdaeq}
\int_{\R^n}|f(x)-g(x+w)|\,dx\leq \omega_\lambda(\varepsilon)\int_{\R^n} f
\end{equation}
where $\omega_\lambda(\varepsilon)=
c^nn^n\left(\frac{\varepsilon}{\min\{\lambda,1-\lambda\}}\right)^{\frac1{19}}$
for some absolute constant $c>1$.
\end{theo}

We frequently use the notation 
$\R_{\geq 0}=\{x\in\R:\, x\geq 0\}$.

\begin{theo}
\label{coordinatewisestab0}
If $\lambda\in (0,1)$ and unconditional convex bodies $K$ and $C$ in $\R^n$
satisfy
$$
V(K^{1-\lambda}\cdot C^\lambda)\leq (1+\varepsilon) V(K)^{1-\lambda} V(C)^\lambda
$$
for $\varepsilon>0$, then
there exists  positive definite diagonal matrix $\Phi$ such that
\begin{equation}
\label{coordinatewisestab0eq}
V(K\Delta (\Phi C))<8\omega_\lambda(\varepsilon) V(K)
\mbox{ \ and }
V((\Phi^{-1}K)\Delta C)<12\omega_\lambda(\varepsilon) V(C)
\end{equation}
where $\omega_\lambda(\varepsilon)$ is taken from \eqref{PLhstablambdaeq}.
\end{theo}
\proof To simplify notation, for any unconditional convex body $L$, we write
$$
L_+=L\cap \R^n_{\geq 0}.
$$
We may assume that
$$
V(K)=V(C)=1.
$$
If $\omega_\lambda(\varepsilon) \geq \frac14$, then we may choose
$\Phi$ to be any linear map with $\det\Phi=1$, and
$V(K\Delta (\Phi C))< 2$ implies \eqref{coordinatewisestab0eq}. Therefore,
 we may also assume that $\varepsilon>0$ is small enough to ensure
\begin{equation}
\label{omegalambda14}
\omega_\lambda(\varepsilon) <\frac14.
\end{equation}

We set $M=K^{1-\lambda}\cdot C^\lambda$, and consider the log-concave functions
$f,g,h:\,\R^n\to[0,\infty)$ defined by
$f(x_1,\ldots,x_n)=\mathbf{1}_{K_+}(e^{x_1},\ldots,e^{x_n})e^{x_1+\ldots +x_n}$, 
$g(x_1,\ldots,x_n)= \mathbf{1}_{C_+}(e^{x_1},\ldots,e^{x_n})e^{x_1+\ldots +x_n}$ and
$h(x_1,\ldots,x_n)=\mathbf{1}_{M_+}(e^{x_1},\ldots,e^{x_n})e^{x_1+\ldots +x_n}$.
In particular,
$$
h(z)=\sup_{z=(1-\lambda)x+\lambda y} f(x)^{1-\lambda}g(y)^\lambda
$$
holds for any
$z\in\mathbb{R}^n$ by the definition of the coordinatewise product.  In addition, we have
$\int_{\R^n}f=\int_{\R^n}g=\frac1{2^n}$ and
\begin{eqnarray*}
\int_{\R^n}h&=&V(M_+)=\frac{V(M)}{2^n}\leq (1+\varepsilon)
\left(\frac{V(K)}{2^n}\right)^{1-\lambda}\left(\frac{V(C)}{2^n}\right)^\lambda\\
&=& (1+\varepsilon) \left(\int_{\R^n}f\right)^{1-\lambda}\left( \int_{\R^n}g\right)^\lambda=(1+\varepsilon)\int_{\R^n}f.
\end{eqnarray*}
Therefore, Theorem~\ref{PLhstablambda} yields that
 there exists $w=(w_1,\ldots,w_n)\in\R^n$ such that
$$
\int_{\R^n}|f(x)-g(x+w)|\,dx\leq \frac{\omega_\lambda(\varepsilon)}{2^n}=\omega_\lambda(\varepsilon) V(K_+) .
$$
Let $\Phi\in{\rm GL}(n)$ be the diagonal transformation
$\Phi(t_1,\ldots,t_n)=(e^{-w_1}t_1,\ldots,e^{-w_n}t_n)$; therefore,
$$
g(x+w)=a\mathbf{1}_{(\Phi C)_+}(e^{x_1},\ldots,e^{x_n})e^{x_1+\ldots +x_n}=a\tilde{g}(x)
\mbox{ \ where }a=e^{w_1+\ldots +w_n}.
$$
We deduce that
\begin{eqnarray*}
\omega_\lambda(\varepsilon) V(K_+)&\geq &\int_{\R^n}|f(x)-a\tilde{g}(x)|\,dx=
\int_{\R^n_+}|\mathbf{1}_{K_+}-a\mathbf{1}_{(\Phi C)_+}|\\
&=&|a-1|V(K_+\cap (\Phi C)_+)+V(K_+\backslash (\Phi C)_+)+aV((\Phi C)_+\backslash K_+).
\end{eqnarray*}
In particular, we have
\begin{equation}
\label{KminusTC}
V(K_+\backslash (\Phi C)_+)\leq \omega_\lambda(\varepsilon) V(K_+),
\end{equation}
and hence  (\ref{omegalambda14}) implies  that
$V(K_+\cap (\Phi C)_+)\geq \frac34\, V(K_+)$. In turn, we deduce
$$
|a-1|\leq \frac{\omega_\lambda(\varepsilon) V(K_+)}{V(K_+\cap (\Phi C)_+)}\leq
\frac43\, \omega_\lambda(\varepsilon)<\frac13,
$$
thus $a>\frac23$. It follows that
\begin{equation}
\label{TCminusK}
V((\Phi C)_+\backslash K_+)\leq \frac{\omega_\lambda(\varepsilon) V(K_+)}{a}
<\frac32\,\omega_\lambda(\varepsilon) V(K_+).
\end{equation}
Combining (\ref{KminusTC}) and (\ref{TCminusK}) yields
$V(K_+\Delta (\Phi C)_+)<3\omega_\lambda(\varepsilon)V(K_+)$, and hence
$V(K\Delta (\Phi C))<3\omega_\lambda(\varepsilon) V(K)$.

Finally, $V(K\Delta (\Phi C))<3\omega_\lambda(\varepsilon) V(K)$ and $\omega_\lambda(\varepsilon)\leq \frac14$
yield that $V(\Phi C)\geq \frac14\,V(K)$, and hence $V(K\Delta (\Phi C))<12\omega_\lambda(\varepsilon) V(\Phi C)$.
\proofbox

\section{Linear images of unconditional convex bodies}
\label{seclinimage}

The main additional tool in this section is to strengthen the containment relation
$$
 K^{1-\lambda}\cdot C^\lambda\subset (1-\lambda)\cdot K +_0 \lambda\cdot C.
$$
Let $e_1,\ldots,e_n$ be the standard orthonormal basis of $\R^n$ that is used in the definition of unconditionality. For a proper subset $J\subset \{1,\ldots,n\}$,
we set
$$
L_J={\rm lin}\{e_i\}_{i\in J}.
$$
We observe that for a diagonal matrix $T={\rm diag}(t_1,\ldots,t_n)$, its operator norm is
\begin{equation}
\label{DiagOperatorNorm}
\|T\|_\infty=\max_{i=1,\ldots,n}|t_i|.
\end{equation}
We write $B^n$ to denote the unit Euclidean ball centered at the origin,  and $e_1,\ldots,e_n$ to denote the fixed orthonormal basis  of $\R^n$ that we use in the definition of unconditionality. We write
${\rm conv} X$ to denote the convex hull of an $X\subset \R^n$, and
frequently use that if $C\subset\R^n$ is unconditional, then $\pm  h_C(e_i)e_i\in \partial C$
\begin{align}
\label{uncond-sandwich}
{\rm conv}\{\pm h_C(e_i)e_i:\,i=1,\ldots,n\}\subset& C\subset \bigoplus_{i=1}^nh_C(e_i)\cdot {\rm conv}\{-e_i,e_i\}\\
\nonumber
\subset&
n\cdot {\rm conv}\{\pm h_C(e_i)e_i:\,i=1,\ldots,n\}.
\end{align}

\begin{prop}
\label{coordlog}
If $\tau\in(0,\frac12]$, $\lambda\in(\tau,1-\tau)$, $K$ is an unconditional convex body in $\R^n$ and $\Phi$ is a positive definite diagonal matrix
 satisfying
$$
V((1-\lambda)\cdot K +_0 \lambda\cdot (\Phi K))\leq (1+\varepsilon)V(K^{1-\lambda}\cdot (\Phi K)^\lambda)
$$
for $\varepsilon>0$,
then either
$\|s\Phi-I_n\|_\infty\leq 16n^4\cdot \frac{\varepsilon^{\frac1{5n}}}{\tau^{\frac15}}$ for $s=(\det \Phi)^{\frac{-1}n}$, or there exist
$s_1,\ldots,s_m>0$  and a partition
of $\{1,\ldots,n\}$ into proper subsets $J_1,\ldots,J_m$, $m\geq 2$, such that
$$
\bigoplus_{k=1}^m (L_{J_k}\cap K)\subset
\left(1+16n^4\cdot \frac{\varepsilon^{\frac1{5n}}}{\tau^{\frac15}} \right)K
$$
for $L_{J_k}={\rm lin}\{e_i:\, i\in J_k\}$, $k=1,\ldots,m$, and in addition,  we have
$$
s_k\cdot (L_{J_k}\cap K)
\subset \Phi(L_{J_k}\cap K)\subset \left(1+16n^4\cdot \frac{\varepsilon^{\frac1{5n}}}{\tau^{\frac15}} \right)s_k\cdot (L_{J_k}\cap K).
$$
\end{prop}
\proof Let 
 $\Phi={\rm diag}(\alpha_1,\ldots,\alpha_n)$, $\alpha_i>0$.

First we assume that
$\varepsilon\geq\frac{\tau^n}{2^{20n}n^{15n}}$, and hence
$$
16n^4\cdot \frac{\varepsilon^{\frac1{5n}}}{\tau^{\frac15}}\geq n.
$$
We set $m=n$, $J_i=\{i\}$  and $s_i=\alpha_i$ for $i=1,\ldots,n$, thus \eqref{uncond-sandwich}
implies Proposition~\ref{coordlog}.

Therefore, from now on, we assume that
\begin{equation}
\label{epscond}
\varepsilon<\frac{\tau^n}{2^{20n}n^{15n}}.
\end{equation}

Since we may apply a positive definite diagonal transform to $K$, we may also assume that
$$
e_i\in\partial \Phi^\lambda K=\partial(K^{1-\lambda}\cdot (\Phi K)^\lambda)\mbox{ \ for $i=1,\ldots,n$}.
$$

Let
$$
\theta= 8n^2\cdot\frac{\varepsilon^{\frac1{5n}}}{\tau^{\frac15}}<\frac1{2n}.
$$

We write $i\bowtie j$ for $i,j\in\{1,\ldots,n\}$
if
$$
\exp(-\theta)\leq \frac{\alpha_i}{\alpha_j}\leq
\exp(\theta).
$$
In addition, we write $\sim$ to denote the equivalence relation on $\{1,\ldots,n\}$ induced by $\bowtie$; namely,
for $i,j\in \{1,\ldots,n\}$, we have $i\sim j$ if and only if there exist pairwise different
$i_0,\ldots,i_l\in \{1,\ldots,n\}$ with $i_0=i$, $i_l=j$, and $i_{k-1}\bowtie i_k$ for $k=1,\ldots,l$.
We may readily assume that
\begin{equation}
\label{ijdifferent}
l\leq n\mbox{ \ in the definition of $i\sim j$}.
\end{equation}

Let $J_1,\ldots,J_m$, $m\geq 1$ be the equivalence classes with respect to $\sim$.
The reason behind introducing $\sim$ are the estimates \eqref{PhibetakLJ}, (i) and (ii).
We claim that if $k=1,\ldots,m$ and $s_k=\min\{\alpha_i:\,i\in J_k\}$, then any $x\in L_{J_k}$ satisfies
\begin{equation}
\label{PhibetakLJ}
s_k\|x\|\leq \|\Phi x\|\leq
e^{n\,\theta}s_k\|x\|.
\end{equation}
To prove \eqref{PhibetakLJ}, we choose $\tilde{i},\tilde{j}\in J_k$ satisfying
$\alpha_{\tilde{i}}=\min\{\alpha_i:\,i\in J_k\}=s_k$
and $\alpha_{\tilde{j}}=\max\{\alpha_i:\,i\in J_k\}$.
We deduce from \eqref{ijdifferent} that
$\alpha_{\tilde{j}}/\alpha_{\tilde{i}}\leq e^{n\,\theta}$,
and hence $s_k\leq \alpha_i\leq e^{n\,\theta}s_k$
holds for $i\in J_k$, proving \eqref{PhibetakLJ}.

Next, if $k\neq l$ holds for $k,l\in\{1,\ldots,m\}$, then
the definition of the relation $\sim$ yields that
 either
$\min\{\alpha_i:\,i\in J_k\}\geq e^{\theta}\cdot
\max\{\alpha_j:\,i\in J_l\}$, or
$\max\{\alpha_i:\,i\in J_k\}\leq e^{-\theta}\cdot
\min\{\alpha_j:\,i\in J_l\}$; therefore,
\begin{description}
\item[(i)] either
$\frac{\|\Phi x\|}{\|x\|}\geq e^{ \theta}\cdot
\frac{\|\Phi y\|}{\|y\|}$
for any
$x\in L_{J_k}\backslash o$ and
$y\in L_{J_l}\backslash o$;
\item[(ii)]  or
$\frac{\|\Phi x\|}{\|x\|}\leq e^{- \theta}\cdot
\frac{\|\Phi y\|}{\|y\|}$
for any
$x\in L_{J_k}\backslash o$ and
$y\in L_{J_l}\backslash o$.
\end{description}

First we assume that
\begin{equation}
\label{m=1}
m=1.
\end{equation}
For our fixed orthonormal basis $e_1,\ldots,e_n$ of $\R^n$, we have $\Phi e_i=\alpha_ie_i$ for $i=1,\ldots,n$. If $m=1$, and hence $J_1=\{1,\ldots,n\}$, then
\eqref{PhibetakLJ} yields that
$$
e^{-n\theta}s_1^{-1}\leq s=(\det \Phi)^{\frac{-1}n}=\left(\prod_{i=1}^n\alpha_i\right)^{\frac{-1}n}\leq s_1^{-1}.
$$
These estimates together with $n\theta\leq 1$ imply that $1-n\theta\leq s\alpha_i\leq 1+2n\theta$ for $i=1,\ldots,n$, and hence ({\it cf.} \eqref{DiagOperatorNorm})
\begin{equation}
\label{PhiI2ntheta}
\|s\Phi-I_n\|_\infty\leq 2n\theta,
\end{equation}
completing the proof of Proposition~\ref{coordlog} if $m=1$.\\

Therefore, we assume that $m\geq 2$. Here again \eqref{PhibetakLJ} yields that
if $k=1,\ldots,m$, then
\begin{equation}
\label{Phibetaktheta}
s_k\cdot (L_{J_k}\cap K)
\subset \Phi(L_{J_k}\cap K)\subset
(1+2n\theta)s_k\cdot (L_{J_k}\cap K).
\end{equation}
For
\begin{equation}
\label{mdef}
M=\bigoplus_{k=1}^m (L_{J_k}\cap \Phi^\lambda K)\supset  \Phi^\lambda K,
\end{equation}
the condition $e_i\in\partial \Phi^\lambda K$, $i=1,\ldots,n$, the convexity and the unconditionality of $\Phi^\lambda K$ and $M$ yield that
${\rm conv}\{\pm e_1,\ldots,\pm e_n\}\subset M\subset [-1,1]^n$, and hence
\begin{equation}
\label{MBn}
\frac1{\sqrt{n}}\,B^n\subset
M\subset  \sqrt{n}\,B^n.
\end{equation}

The core statement is that
\begin{equation}
\label{MinKpos}
(1-2\sqrt{n}\theta )M\subset \Phi^\lambda K,
\end{equation}
what would complete the proof of Proposition~\ref{coordlog} by \eqref{mdef}.

We prove \eqref{MinKpos} indirectly, and hence we suppose that
\begin{equation}
\label{MnotinK}
(1-2\sqrt{n}\theta )M\not\subset
\Phi^\lambda K,
\end{equation}
and seek a contradiction in three steps.\\

\noindent{\bf Step 1} {\it If $m\geq 2$ and the indirect hypothesis \eqref{MnotinK} holds,
then we indentify an $x_0\in  \partial(\Phi^\lambda K)\cap \R^n_{\geq 0}$ in Step~1, which
sits "deep" in $(1-\lambda)\cdot K+_0 \lambda\cdot (\Phi K)$ according to Claim~\ref{x0neighborhoodclaim} in Step~2.}

Let $\eta>0$ be maximal such that
$$
\eta(M+\theta B^n)\subset \Phi^\lambda K.
$$
We deduce that
\begin{equation}
\label{etabounds}
\frac1{2n}\leq \eta <1-2\sqrt{n}\theta.
\end{equation}
where the upper bound follows from \eqref{MnotinK}, and the lower bound
follows from $\frac1n M\subset \Phi^\lambda K$ (as $ \Phi^\lambda K$ unconditional)
and 
$\theta B^n\subset M$ that is the consequence of \eqref{MBn}.

The maximality of $\eta$ and the unconditionality of $K$ yield that there exists an
$$
x_0\in \eta(M+\theta B^n)\cap  \partial(\Phi^\lambda K)\cap \R^n_{\geq 0},
$$
and there exists a unique exterior normal $w\in S^{n-1}\cap \R^n_{\geq 0}$
to $\partial(\Phi^\lambda K)$ at $x_0$ satisfying
({\it cf.} \eqref{etabounds})
\begin{equation}
\label{x0ball}
x_0-\frac{\theta}{2n}\cdot w+\frac{\theta}{2n}\cdot B^n\subset \Phi^\lambda K.
\end{equation}
In addition, we have
\begin{equation}
\label{x0plusball}
x_0+\theta\,B^n\subset \eta(M+\theta B^n)+\theta B^n\subset (1-2\sqrt{n}\theta)M+2\theta B^n
\subset M.
\end{equation}

Writing $x|L$ to denote the orthogonal projection of $x\in\R^n$ to a linear subspace $L$, we claim that
\begin{equation}
\label{wprojup}
\|w|L_{J_k}\|^2\leq 1-\frac{\theta^2}{2n}\mbox{ \ for $k=1,\ldots,m$}.
\end{equation}
Let $v\in S^{n-1}\cap L_{J_k}$
satisfy $w|L_{J_k}=\|w|L_{J_k}\|\,v$, and hence
$$
\|w|L_{J_k}\|=\langle w,v\rangle.
$$
Since $\|x_0\|\leq \sqrt{n}$ by \eqref{MBn} and
$x_0-(x_0|L_{J_k})$ is orthogonal to $v$, we have
$$
|\langle w,x_0-(x_0|L_{J_k})\rangle|=\|x_0-(x_0|L_{J_k})\|\sqrt{1-\langle w,v\rangle^2}\leq
\sqrt{n}\sqrt{1-\langle w,v\rangle^2}.
$$
It follows from \eqref{x0plusball} that
$$
(x_0|L_{J_k})+\theta v\in \Phi^\lambda K\cap L_{J_k}.
$$
Since $w$ is an exterior normal to $\Phi^\lambda  K$ at $x_0$, we have
$\langle w,x_0\rangle\geq \langle w,(x_0|L_{J_k})+\theta v\rangle$, thus
$$
\sqrt{n}\sqrt{1-\langle w,v\rangle^2}\geq \langle w,x_0-(x_0|L_{J_k})\rangle
\geq \theta\langle w, v\rangle.
$$
We deduce that
$$
\|w|L_{J_k}\|^2=\langle w,v\rangle^2\leq \frac{n}{n+\theta^2}=1-\frac{\theta^2}{n+\theta^2}<
1-\frac{\theta^2}{2n},
$$
proving \eqref{wprojup}.

In turn, we conclude from
$\sum_{k=1}^m\|w|L_{J_k}\|^2=1$,
$m\leq n$ and \eqref{wprojup} that
there exist $p\neq q$ satisfying
$$
\|w|L_{J_p}\|^2\geq \frac{\theta^2}{2n^2}\mbox{ and }
\|w|L_{J_q}\|^2\geq \frac{\theta^2}{2n^2}.
$$
Possibly after reindexing, we may assume that
\begin{equation}
\label{wprojLJlow}
\|w|L_{J_1}\|\geq \frac{\theta}{2n}\mbox{ and }
\|w|L_{J_2}\|\geq \frac{\theta}{2n}.
\end{equation}

For any $u\in S^{n-1}\cap \R^n_{\geq 0}$, it follows from $\Phi^{-\lambda}x_0\in K$ and $\Phi^{1-\lambda}x_0\in \Phi K$  that
\begin{equation}
\label{x0neighborhood01}
\langle u, \Phi^{-\lambda}x_0\rangle\leq h_K(u)\mbox{ \ and \ }
\langle u, \Phi^{1-\lambda}x_0\rangle\leq h_{\Phi K}(u);
\end{equation}
therefore the H\"older inequality yields that 
\begin{equation}
\label{x0neighborhood0}
\langle u,x_0\rangle\leq \langle u, \Phi^{-\lambda}x_0\rangle^{1-\lambda}\langle u, \Phi^{1-\lambda}x_0\rangle^{\lambda}
\leq h_K(u)^{1-\lambda}h_{\Phi K}(u)^{\lambda}.
\end{equation}
In particular,  \eqref{x0neighborhood0} implies that $x_0\in (1-\lambda)\cdot K+_0 \lambda\cdot (\Phi K)$.

In order to prove \eqref{MinKpos}; more precisely, to prove that
\eqref{MnotinK} is false, the next step is the following stability version of \eqref{x0neighborhood0}.\\

\noindent{\bf Step 2} {\it If $m\geq 2$ and the indirect hypothesis \eqref{MnotinK} holds, then the $x_0$ from Step~1
sits "deep" in $(1-\lambda)\cdot K+_0 \lambda\cdot (\Phi K)$.}

\begin{claim}
\label{x0neighborhoodclaim}
For any $u\in S^{n-1}\cap \R^n_{\geq 0}$, we have
\begin{equation}
\label{x0neighborhood}
\langle u,x_0\rangle\left(1+\frac{\tau\theta^5}{1024n^{5.5}}\right)\leq h_K(u)^{1-\lambda}h_{\Phi K}(u)^{\lambda}.
\end{equation}
\end{claim}
\proof We observe that
$\langle u, \Phi^{-\lambda}x_0\rangle=\langle \Phi^{-\lambda}u, x_0\rangle$,
 $\langle u, \Phi^{1-\lambda}x_0\rangle=\langle \Phi^{1-\lambda}u, x_0\rangle$,
\begin{eqnarray*}
 h_K(u)&=&h_{\Phi^{\lambda} K}(\Phi^{-\lambda}u);\\
 h_{\Phi K}(u)&=&h_{\Phi^{\lambda} K}(\Phi^{1-\lambda}u),
\end{eqnarray*}
and hence it follows from \eqref{x0neighborhood01} and
 \eqref{x0neighborhood0} that it is sufficient to prove that if $u\in S^{n-1}$, then either
\begin{equation}
\label{x0neighborhoodsuff}
\mbox{$\left(\frac{h_{\Phi^{\lambda} K}(\Phi^{-\lambda}u)}{\langle \Phi^{-\lambda} u,x_0\rangle}\right)^{1-\lambda}
\geq 1+\frac{\tau\theta^5}{1024n^{5.5}}$,
 or
$\left(\frac{h_{\Phi^{\lambda} K}(\Phi^{1-\lambda}u)}{\langle \Phi^{1-\lambda} u,x_0\rangle}\right)^{\lambda}\geq 1+\frac{\tau\theta^5}{1024n^{5.5}}$.}
\end{equation}
Let us write $w=\bigoplus_{k=1}^mw_k$ and $u=\bigoplus_{k=1}^m u_k$
for $w_k=w|L_{J_k}$ and $u_k=u|L_{J_k}$, and prove that ({\it cf.} \eqref{wprojLJlow})
there exists $i\in\{1,2\}$ such that
\begin{equation}
\label{x0neighborhoodw12}
\mbox{either \ }
\left|\frac{\|\Phi^{-\lambda}u_i\|}{\|\Phi^{-\lambda}u\|}-\|w_i\|\right|\geq \frac{\theta^2}{16n^2},
\mbox{ \ or \ } \left|\frac{\|\Phi^{1-\lambda}u_i\|}{\|\Phi^{1-\lambda}u\|}-\|w_i\|\right|\geq \frac{\theta^2}{16n^2}.
\end{equation}
We prove
\eqref{x0neighborhoodw12}  by contradiction; thus, we suppose that
if $i\in\{1,2\}$, then
$$
\left|\frac{\|\Phi^{-\lambda}u_i\|}{\|\Phi^{-\lambda}u\|}-\|w_i\|\right|< \frac{\theta^2}{16n^2}
\mbox{ \ and \ } \left|\frac{\|\Phi^{1-\lambda}u_i\|}{\|\Phi^{1-\lambda}u\|}-\|w_i\|\right|< \frac{\theta^2}{16n^2}.
$$
and seek a contradiction. Since $\|w_1\|\geq \frac{\theta}{2n}$ and
$\|w_2\|\geq \frac{\theta}{2n}$ according to \eqref{wprojLJlow}, we deduce that
if $i\in\{1,2\}$, then
\begin{equation}
\label{x0neighborhoodindirect12}
e^{-\frac{\theta}4}<\frac{\|\Phi^{-\lambda}u_i\|}{\|\Phi^{-\lambda}u\|\cdot \|w_i\|}< e^{\frac{\theta}4},
\mbox{ \ and \ }e^{-\frac{\theta}4}<\frac{\|\Phi^{1-\lambda}u_i\|}{\|\Phi^{1-\lambda}u\|\cdot \|w_i\|}< e^{\frac{\theta}4}.
\end{equation}
It follows from $\Phi\left(\Phi^{-\lambda}u_1\right)=\Phi^{1-\lambda}u_1$,
$\Phi\left(\Phi^{-\lambda}u_2\right)=\Phi^{1-\lambda}u_2$,
 and \eqref{x0neighborhoodindirect12} that
$$
e^{-\frac{\theta}2}<\frac{\|\Phi\left(\Phi^{-\lambda}u_1\right)\|}{\|\Phi^{-\lambda}u_1\|}< e^{\frac{\theta}2}
\mbox{ \ and \ }
e^{-\frac{\theta}2}<\frac{\|\Phi\left(\Phi^{-\lambda}u_2\right)\|}{\|\Phi^{-\lambda}u_2\|}< e^{\frac{\theta}2};
$$
therefore,
$$
e^{-\theta}<\frac{\|\Phi\left(\Phi^{-\lambda}u_1\right)\|}{\|\Phi^{-\lambda}u_1\|}:
\frac{\|\Phi\left(\Phi^{-\lambda}u_2\right)\|}{\|\Phi^{-\lambda}u_2\|}< e^{\theta}.
$$
Since $\Phi^{-\lambda}u_i\in L_{J_i}$ for $i=1,2$, the last inequalities contradict (i) and (ii),
and in turn verify \eqref{x0neighborhoodw12}.

Based on \eqref{x0neighborhoodw12}, the triangle inequality yields
the existence of $i\in\{1,2\}$ such that
$$
\mbox{either \ }
\left\|\frac{\Phi^{-\lambda}u_i}{\|\Phi^{-\lambda}u\|}-w_i\right\|\geq \frac{\theta^2}{16n^2}
\mbox{ \ or \ }
\left\|\frac{\Phi^{1-\lambda}u_i}{\|\Phi^{1-\lambda}u\|}-w_i\right\|\geq \frac{\theta^2}{16n^2},
$$
and in turn we deduce that
\begin{equation}
\label{Philambda1u}
\mbox{either \ }\left\|\frac{\Phi^{-\lambda}u}{\|\Phi^{-\lambda}u\|}-w\right\|\geq \frac{\theta^2}{16n^2}
\mbox{ \ or \ }
\left\|\frac{\Phi^{1-\lambda}u}{\|\Phi^{1-\lambda}u\|}-w\right\|\geq \frac{\theta^2}{16n^2}.
\end{equation}

First, we assume that out of the two possibilities in \eqref{Philambda1u}, we have
\begin{equation}
\label{Philambdau}
\left\|\frac{\Phi^{-\lambda}u}{\|\Phi^{-\lambda}u\|}-w\right\|\geq \frac{\theta^2}{16n^2}.
\end{equation}
According to \eqref{x0ball}, we have
$$
\widetilde{B}=x_0-\frac{\theta}{2n}\cdot w+\frac{\theta}{2n}\cdot B^n\subset \Phi^\lambda K,
$$
which in turn yields (using  that $\langle v,v-w\rangle=\frac12\,\|v-w\|^2$ for $v\in S^{n-1}$ as $w\in S^{n-1}$,
and later  \eqref{Philambdau} and $\|x_0\|\leq\sqrt{n}$) that
\begin{align*}
h_{\Phi^{\lambda}K}(\Phi^{-\lambda}u)-\langle \Phi^{-\lambda}u,x_0\rangle\geq&
h_{\widetilde{B}}(\Phi^{-\lambda}u)-\langle \Phi^{-\lambda}u,x_0\rangle\\
=& \left\langle \Phi^{-\lambda}u, x_0-\frac{\theta}{2n}\cdot w+
\frac{\theta}{2n}\cdot \frac{\Phi^{-\lambda}u}{\|\Phi^{-\lambda}u\|}\right\rangle
-\langle \Phi^{-\lambda}u,x_0\rangle\\
=& \|\Phi^{-\lambda}u\|\cdot \frac{\theta}{4n}\cdot
\left\|\frac{\Phi^{-\lambda}u}{\|\Phi^{-\lambda}u\|}-w\right\|^2\\
\geq & \frac{\langle \Phi^{-\lambda}u,x_0\rangle}{\sqrt{n}}\cdot \frac{\theta}{4n}
\left(\frac{\theta^2}{16n^2}\right)^2= \frac{\langle \Phi^{-\lambda}u,x_0\rangle\cdot\theta^5}{1024n^{5.5}}.
\end{align*}
We conclude using $1-\lambda\geq \tau$ that
\begin{equation}
\label{x0neighborhoodfirst}
\left(\frac{h_{\Phi^{\lambda} K}(\Phi^{-\lambda}u)}{\langle \Phi^{-\lambda} u,x_0\rangle}\right)^{1-\lambda}\geq
\left(\frac{h_{\Phi^{\lambda} K}(\Phi^{-\lambda}u)}{\langle \Phi^{-\lambda} u,x_0\rangle}\right)^{\tau}
\geq 1+\frac{\tau\theta^5}{1024n^{5.5}}.
\end{equation}

Secondly, if
$$
\left\|\frac{\Phi^{1-\lambda}u}{\|\Phi^{1-\lambda}u\|}-w\right\|\geq \frac{\theta^2}{16n^2}
$$
holds in \eqref{Philambda1u}, then similar argument yields
$$
\left(\frac{h_{\Phi^{\lambda} K}(\Phi^{1-\lambda}u)}{\langle \Phi^{1-\lambda} u,x_0\rangle}\right)^{\lambda}
\geq 1+\frac{\tau\theta^5}{1024n^{5.5}}.
$$
proving \eqref{x0neighborhoodsuff}. In turn, we conclude \eqref{x0neighborhood}
in Claim~\ref{x0neighborhoodclaim}.
\proofbox

\noindent{\bf Step 3} {\it If $m\geq 2$ and the indirect hypothesis \eqref{MnotinK} holds, then Claim~\ref{x0neighborhoodclaim} contradicts  the assumption} 
$$
V((1-\lambda)\cdot K +_0 \lambda\cdot (\Phi K))\leq (1+\varepsilon)V(K^{1-\lambda}\cdot (\Phi K)^\lambda).
$$

Let $\varrho\geq 0$ be maximal with the property that
\begin{equation}
\label{rhodef}
x_0+\varrho B^n\subset (1-\lambda)\cdot K+_0 \lambda\cdot (\Phi K).
\end{equation}
We claim that
\begin{equation}
\label{rholow}
\varrho\geq \frac{\tau\theta^5}{2048n^{6}}.
\end{equation}
It follows from Claim~\ref{x0neighborhoodclaim} that $\varrho>0$. To prove \eqref{rholow}, we may assume that
\begin{equation}
\label{rhoupass}
\varrho\leq \frac{\tau\theta^5}{2048n^{6}}<\frac1{2n}.
\end{equation}
We consider a
$$
y_0\in (x_0+\varrho B^n)\cap\partial \big( (1-\lambda)\cdot K+_0 \lambda\cdot (\Phi K)\big)\cap \R_{\geq 0}^n,
$$
which exists as $(1-\lambda)\cdot K+_0 \lambda\cdot (\Phi K)$ is unconditional. It follows from \eqref{rhodef} that $y_0$ is a smooth boundry point with a unique exterior unit normal 
$u\in S^{n-1}\cap \R_{\geq 0}^n$  to
$$
\widetilde{M}=(1-\lambda)\cdot K+_0 \lambda\cdot (\Phi K)
$$
at $y_0$, and hence $y_0=x_0+\varrho\,u$. On the one hand,
$\pm e_i\in\widetilde{M}$ for $i=1,\ldots,n$ yields that
$h_{\widetilde{M}}(u)\geq \frac1{\sqrt{n}}$, thus \eqref{rhoupass} implies
\begin{equation}
\label{ux0low}
\langle u,x_0\rangle=\langle u,y_0\rangle-\varrho=h_{\widetilde{M}}(u)-\varrho\geq \frac1{2\sqrt{n}}.
\end{equation}
On the other hand, $h_{\widetilde{M}}(u)=h_K(u)^{1-\lambda}h_{\Phi K}(u)^{\lambda}$ holds because $y_0$ is a smooth boundary point of $\widetilde{M}$; therefore, we deduce from Claim~\ref{x0neighborhoodclaim},  \eqref{rhodef}, \eqref{rhoupass} and
\eqref{ux0low} that
\begin{eqnarray*}
\varrho&=&h_{\widetilde{M}}(u)-\langle u,x_0\rangle=
h_K(u)^{1-\lambda}h_{AK}(u)^{\lambda}-\langle u,x_0\rangle\\
&\geq & \langle u,x_0\rangle\cdot \frac{\tau\theta^5}{1024n^{5.5}}\geq \frac{\tau\theta^5}{2048n^{6}},
\end{eqnarray*}
proving \eqref{rholow}.

Since $V(\Phi^\lambda K)\leq 2^n$ because of $\pm e_i\in (\Phi^\lambda K)$, $i=1,\ldots,n$,
$\kappa_n=\frac{\pi^{\frac{n}2}}{\Gamma(\frac{n}2+1)}>\frac{(\pi\,e)^{\frac{n}2}}{4\sqrt{n}\cdot n^{\frac{n}2}}$,
and  the supporting hyperplane at $x_0$ to $\Phi^\lambda K$ cuts $x_0+\varrho B^n$ into half, we deduce that
\begin{eqnarray*}
 V(\widetilde{M})&\geq&
V(\Phi^\lambda K)+\frac{\varrho^n\kappa_n}2\geq
V(\Phi^\lambda K)+\frac{\kappa_n\tau^n\theta^{5n}}{2\cdot 2048^nn^{6n}}\\
&=&
V(\Phi^\lambda K)\left(1+\frac{\kappa_n\tau^n\theta^{5n}}{2\cdot 2048^nn^{6n}V(\Phi^\lambda K)}\right)\\
&>&V(\Phi^\lambda K)\left(1+\frac{(\pi e)^{\frac{n}2}\tau^n\theta^{5n}}{8\sqrt{n}\cdot 4096^nn^{6.5n}}\right)
>V(\Phi^\lambda K)\left(1+\frac{\tau^n\theta^{5n}}{2^{15n}n^{10n}}\right)\\
&>&(1+\varepsilon)V(\Phi^\lambda K)=(1+\varepsilon)V(K)^{1-\lambda}V(\Phi K)^\lambda,
\end{eqnarray*}
what is absurd. This contradicts \eqref{MnotinK}, and  completes Step~3. In particular, this contradiction verifies
$(1-2\sqrt{n}\theta )M\subset \Phi^\lambda K$ ({\it cf.} \eqref{MinKpos}), and hence finally
proves Proposition~\ref{coordlog}. 
\proofbox

\section{Proof of Theorem~\ref{logBMuncondstab}}
\label{seclogBMuncondstab}

The proof of Theoem~\ref{logBMuncondstab} will be based on Theorem~\ref{coordinatewisestab0}
and Proposition~\ref{coordlog}. However, first we need some simple lemmas.
The first statement is the following corollary (see Lemma~3.1 of Kolesnikov, Milman \cite{KoM22}) of
the logarithmic Brunn-Minowski inequality for unconditional convex bodies due to Saroglou \cite{Sar15}.

\begin{lemma}
\label{logsum-logconcave}
If $K$ and $C$ are unconditional convex bodies in $\R^n$, then
$$
\varphi(t)=V((1-t)\cdot K +_0 t\cdot C)
$$
is log-concave on $[0,1]$.
\end{lemma}

The second claim provides simple estimates about log-concave functions.

\begin{lemma}
\label{log-concave-error}
Let $\varphi$ be a log-concave function on $[0,1]$.
\begin{description}
\item[(i)] If $\lambda\in(0,1)$, $\eta\in(0,2\cdot\min\{1-\lambda,\lambda\})$ and
  $\varphi(\lambda)\leq (1+\eta)\varphi(0)^{1-\lambda}\varphi(1)^\lambda$, then
$$
\varphi\left(\mbox{$\frac12$}\right)\leq \left(1+\frac{\eta}{\min\{1-\lambda,\lambda\}}\right)\sqrt{\varphi(0)\varphi(1)}
$$
\item[(ii)] If $\varphi(0)=\varphi(1)=1$ and $\varphi'(0)\leq 2$, then
$\varphi\left(\mbox{$\frac12$}\right)\leq 1+\varphi'(0)$.
\end{description}
\end{lemma}
\proof For (i), we  may assume that $0<\lambda<\frac12$, and hence $\lambda=(1-2\lambda)\cdot 0+2\lambda\cdot\frac12$,
$\varphi(\lambda)\leq (1+\eta)\varphi(0)^{1-\lambda}\varphi(1)^\lambda$
and the log-concavity of $\varphi$ yield
$$
(1+\eta)\varphi(0)^{1-\lambda}\varphi(1)^\lambda\geq \varphi(\lambda)\geq
\varphi(0)^{1-2\lambda}\varphi\left(\mbox{$\frac12$}\right)^{2\lambda}.
$$
Thus $(1+\eta)^{\frac1{2\lambda}}\leq e^{\frac{\eta}{2\lambda}}\leq 1+\frac{\eta}{\lambda}$ implies
$$
\varphi\left(\mbox{$\frac12$}\right)\leq (1+\eta)^{\frac1{2\lambda}}\sqrt{\varphi(0)\varphi(1)}\leq
\left(1+\frac{\eta}{\lambda}\right)\sqrt{\varphi(0)\varphi(1)}.
$$
For (ii), we write $\varphi(t)=e^{W(t)}$ for a concave function $W$ with $W(0)=W(1)=0$. Thus
$W(\frac12)\leq \frac12 W'(0)$, which in turn yields using $W'(0)=\varphi'(0)\leq 2$ that
$$
\varphi\left(\mbox{$\frac12$}\right)=e^{W(\frac12)}\leq e^{W'(0)/2}\leq 1+W'(0)=1+\varphi'(0).
\mbox{ \ }\proofbox
$$

We also need the following statement about volume difference.

\begin{lemma}
\label{voldiffhmothete}
If $M\subset K$ are $o$-symmetric convex bodies with $V(K\backslash M))\leq \frac1{2^{n}}\,V(K)$, then
$$
K\subset \left(1+4\cdot \left(\frac{V(K\backslash M)}{V(M)}\right)^{\frac1n}\right) M.
$$
\end{lemma}
\proof Let $t\geq 0$ be minimal with
$$
K\subset (1+t)M.
$$
Then there exist $z\in \partial K$ and $y\in\partial M$ with $z=(1+t)y$. We have
$$
\frac{2}{2+t}\cdot z=\frac{2(1+t)}{2+t}\cdot y\not \in \frac{2(1+t)}{2+t}\cdot {\rm int}\,M={\rm int}\,M-\frac{t}{2+t}\cdot M
$$
and $\frac{2}{2+t}\cdot z+\frac{t}{2+t}\cdot M\subset K$, and hence
$$
\frac{2}{2+t}\cdot z+\frac{t}{2+t}\cdot M\subset K\backslash {\rm int}\,M.
$$
It follows that
$V(K\backslash M)\geq \left(\frac{t}{2+t}\right)^{n}\cdot V(M)$,
which, together with $V(K\backslash M))\leq \frac1{2^n}\,V(M)$, implies
$t\leq 4\cdot \left(\frac{V(K\backslash M)}{V(M)}\right)^{\frac1n}$.
\proofbox

We will need the case $\lambda=\frac12$ of Theorem~\ref{coordinatewisestab0}
and Proposition~\ref{coordlog}.

\begin{coro}
\label{coordinatewisestab00}
If the unconditional convex bodies $K$ and $C$ in $\R^n$
satisfy
$$
V(K^{\frac12}\cdot C^{\frac12})\leq (1+\varepsilon) V(K)^{\frac12} V(C)^{\frac12}
$$
for $\varepsilon>0$, then
there exists  positive definite diagonal matrix $\Phi$ such that
\begin{equation}
\label{coordinatewisestab00eq}
V(K\Delta (\Phi C))<c^nn^n\varepsilon^{\frac1{19}} V(K)
\end{equation}
where $c>1$ is an absolute constant.
\end{coro}

\begin{coro}
\label{coordlog0}
If $K$ is an unconditional convex body in $\R^n$ and $\Phi$ is a positive definite diagonal matrix
 satisfying
$$
V\left(\frac12\cdot K +_0 \frac12\cdot (\Phi K)\right)\leq (1+\varepsilon)V(K^{\frac12}\cdot (\Phi K)^{\frac12})
$$
for $\varepsilon>0$,
then either
$\|s\Phi-I_n\|_\infty\leq 20n^4\cdot \varepsilon^{\frac1{5n}}$ for $s=(\det \Phi)^{\frac{-1}n}$, or there exist
$s_1,\ldots,s_m>0$  and a partition
of $\{1,\ldots,n\}$ into proper subsets $J_1,\ldots,J_m$, $m\geq 2$, such that
\begin{align*}
\bigoplus_{k=1}^m (L_{J_k}\cap K)&\subset
\left(1+20n^4\cdot \varepsilon^{\frac1{5n}} \right)K\\
s_k(L_{J_k}\cap K)\subset
 \Phi(L_{J_k}\cap K)&\subset
\left(1+20n^4\cdot \varepsilon^{\frac1{5n}}\right)s_k (L_{J_k}\cap K),\;\;
k=1,\ldots,m.
\end{align*}
\end{coro}

\noindent{\bf Proof of Theorem~\ref{logBMuncondstab} }
First we consider the case $\lambda=\frac12$, and hence prove that
if the unconditional convex bodies $K$ and $C$ in $\R^n$
satisfy
\begin{equation}
\label{12cond}
V\left(\frac12\cdot K +_0 \frac12\cdot C\right)\leq (1+\varepsilon) V(K)^{\frac12} V(C)^{\frac12}
\end{equation}
for $\varepsilon>0$, then for $m\geq 1$,
there exist $\theta_1,\ldots,\theta_m>0$ and unconditional compact convex sets
$K_1,\ldots,K_m>0$ such that ${\rm lin}\,K_i$, $i=1,\ldots,m$, are complementary coordinate subspaces, and
\begin{eqnarray}
\label{12condK}
K_1\oplus\ldots \oplus K_m\subset &K&\subset
\left(1+c_0^n\varepsilon^{\frac1{95n}}\right)(K_1\oplus\ldots \oplus K_m)\\
\label{12condC}
\theta_1K_1\oplus\ldots \oplus \theta_mK_m\subset &C&\subset
\left(1+c_0^n\varepsilon^{\frac1{95n}}\right)(\theta_1K_1\oplus\ldots \oplus \theta_mK_m)
\end{eqnarray}
where $c_0>1$ is an absolute constant.

First we assume that
\begin{equation}
\label{epscondlogBMstab}
\varepsilon<\gamma^{-n}n^{-19n}
\end{equation}
for a suitable absolute constant $\gamma>1$ where $\gamma$ is a chosen in a way such that
\begin{equation}
\label{gammacondlogBMstab}
\tilde{c}^nn^n\varepsilon^{\frac1{19}}<\frac1{2^{n}}
\end{equation}
for the constant $\tilde{c}$ of Corollary~\ref{coordinatewisestab00}.

We have
$$
V(K^{\frac12}\cdot C^{\frac12})\leq
V\left(\frac12\cdot K +_0 \frac12\cdot C\right)\leq(1+\varepsilon) V(K)^{\frac12} V(C)^{\frac12};
$$
therefore, Corollary~\ref{coordinatewisestab00} yields a
positive definite diagonal matrix $\Phi$ such that
\begin{equation}
\label{coordinatewisestab0eq1}
V((\Phi K)\Delta C)<\tilde{c}^nn^n\varepsilon^{\frac1{19}} V(C)\mbox{ \ and \ }
V(K\Delta (\Phi^{-1} C))<\tilde{c}^nn^n\varepsilon^{\frac1{19}} V(K)
\end{equation}
where $\tilde{c}>1$ is an absolute constant.

Let
$$
M=K\cap (\Phi^{-1} C),
$$
and hence \eqref{coordinatewisestab0eq1} yields that
\begin{eqnarray}
\label{VMlowK}
V(M)&>&(1- \tilde{c}^nn^n\varepsilon^{\frac1{19}})V(K)\\
\label{VMlowC}
V(\Phi M)&>&(1- \tilde{c}^nn^n\varepsilon^{\frac1{19}})V(C).
\end{eqnarray}
As $M\subset K$ and $\Phi M\subset C$, it follows that
\begin{align*}
\mbox{$V\left(\frac12\, M +_0 \frac12\, (\Phi M)\right)$}\leq&(1+\varepsilon) V(K)^{\frac12} V(C)^{\frac12}\\
\leq& (1+2 \tilde{c}^nn^n\varepsilon^{\frac1{19}}) V(M)^{\frac12} V(\Phi M)^{\frac12}\\
=& (1+2 \tilde{c}^nn^n\varepsilon^{\frac1{19}}) V(M^{\frac12}\cdot(\Phi M)^{\frac12}).
\end{align*}
Now we apply Corollary~\ref{coordlog0}, and deduce the existence of an absolute constant $c_1>0$ such that
either
$\|s\Phi-I_n\|_\infty\leq c_1n^5\cdot \varepsilon^{\frac1{95n}}$ for $s=(\det \Phi)^{\frac{-1}n}$, or there exist
$s_1,\ldots,s_m>0$  and a partition
of $\{1,\ldots,n\}$ into proper subsets $J_1,\ldots,J_m$, $m\geq 2$, such that
$$
\bigoplus_{k=1}^m (L_{J_k}\cap M)\subset
\left(1+c_1n^5\cdot \varepsilon^{\frac1{95n}} \right)M
$$
where  for $k=1,\ldots,m$, we have
$$
s_k\cdot (L_{J_k}\cap M)
\subset \Phi(L_{J_k}\cap M)\subset
\left(1+c_1n^5\cdot \varepsilon^{\frac1{95n}}\right)s_k\cdot (L_{J_k}\cap M).
$$
We deduce from \eqref{gammacondlogBMstab}, \eqref{VMlowK}, \eqref{VMlowC},
and
Lemma~\ref{voldiffhmothete} the existence of an absolute constant $c_2>1$ that
\begin{eqnarray*}
M&\subset K\subset &(1+c_2n\varepsilon^{\frac1{19n}})M\\
\Phi M&\subset C\subset &(1+c_2n\varepsilon^{\frac1{19n}})\Phi M.
\end{eqnarray*}

Now if $\|s\Phi-I_n\|_\infty\leq c_1n^5\cdot \varepsilon^{\frac1{95n}}$, then we can choose
$m=1$ and $K_1=M$ to verify Theorem~\ref{logBMuncondstab}.
On the other hand, if $\|s\Phi-I_n\|_\infty> c_1n^5\cdot \varepsilon^{\frac1{95n}}$, then
we choose
$$
K_k=\left(1+c_1n^5\cdot \varepsilon^{\frac1{95n}} \right)^{-1}(L_{J_k}\cap M)
\mbox{ \ for $k=1,\ldots,m$}.
$$
For $c_3=c_1+c_2+c_1c_2$ and $c_4=c_1+c_3+c_1c_3$,
it follows using $n\varepsilon^{\frac1{19n}}<1$ ({\it cf.} \eqref{epscondlogBMstab}) that
\begin{align*}
\bigoplus_{k=1}^m K_k\subset& M\subset K\subset
(1+c_2n\varepsilon^{\frac1{19n}})M\subset
(1+c_2n\varepsilon^{\frac1{19n}})\bigoplus_{k=1}^m(L_{J_k}\cap M)\\
\subset& \left(1+c_3n^5\cdot \varepsilon^{\frac1{95n}} \right)\bigoplus_{k=1}^m K_k\\
\bigoplus_{k=1}^m s_kK_k\subset&\bigoplus_{k=1}^m \Phi K_k\subset\Phi M\subset C\subset (1+c_2n\varepsilon^{\frac1{19n}})\Phi M\\
\subset& (1+c_2n\varepsilon^{\frac1{19n}})\bigoplus_{k=1}^m\Phi (L_{J_k}\cap M)\\
\subset& (1+c_2n\varepsilon^{\frac1{19n}})\bigoplus_{k=1}^m
\left(1+c_1n^5\cdot \varepsilon^{\frac1{95n}}\right)s_k\cdot (L_{J_k}\cap M)\\
\subset& \left(1+c_3n^5\cdot \varepsilon^{\frac1{95n}}\right)\bigoplus_{k=1}^m
s_k (L_{J_k}\cap M)\subset \left(1+c_4n^5\cdot \varepsilon^{\frac1{95n}}\right)\bigoplus_{k=1}^m
s_k K_k.
\end{align*}
This proves Theorem~\ref{logBMuncondstab} if $\lambda=\frac12$ and
$\varepsilon<\gamma^{-n}n^{-19n}$ ({\it cf.} \eqref{epscondlogBMstab}).

Still keeping  $\lambda=\frac12$, we observe that if $Q$ is any unconditional convex body in $\R^n$, then
\begin{equation}
\label{approximate-box}
\bigoplus_{i=1}^n(\R e_i\cap Q)\subset n Q.
\end{equation}
Therefore, if $\varepsilon\geq \gamma^{-n}n^{-19n}$ ({\it cf.} \eqref{epscondlogBMstab}) holds in \eqref{12cond},
then \eqref{12condK} and \eqref{12condC} readily hold for suitable absolute constant $c_0>1$ by taking
$m=n$, $K_k=\frac1n(\R e_k\cap K)$, and choosing $\theta_k>0$ in a way such that
$\theta_k(\R e_k\cap K)=\R e_k\cap C$ for $k=1,\ldots,n$.
In particular, Theorem~\ref{logBMuncondstab} has been verified if $\lambda=\frac12$.\\

Next, we assume that $\lambda\in[\tau,1-\tau]$ holds for some $\tau\in(0,\frac12]$
in Theorem~\ref{logBMuncondstab}.
First let $\varepsilon\leq \tau$. Since
$$
\varphi(t)=V((1-t)\cdot K +_0 t\cdot C)
$$
is log-concave on $[0,1]$ according to  Lemma~\ref{logsum-logconcave},
Lemma~\ref{log-concave-error} yields that
$$
\varphi\left(\mbox{$\frac12$}\right)\leq \left(1+\frac{\varepsilon}{\min\{1-\lambda,\lambda\}}\right)\sqrt{\varphi(0)\varphi(1)};
$$
or in other words,
$$
V\left(\frac12\cdot K +_0 \frac12\cdot C\right)\leq \left(1+\frac{\varepsilon}{\tau}\right) V(K)^{\frac12} V(C)^{\frac12}.
$$
We deduce from \eqref{12condK} and \eqref{12condC} that
for $m\geq 1$,
there exist $\theta_1,\ldots,\theta_m>0$ and unconditional compact convex sets
$K_1,\ldots,K_m>0$ such that ${\rm lin}\,K_i$, $i=1,\ldots,m$, are complementary coordinate subspaces, and
\begin{eqnarray}
\label{12condKtau}
K_1\oplus\ldots \oplus K_m\subset &K&\subset
\left(1+c_0^n\left(\frac{\varepsilon}{\tau}\right)^{\frac1{95n}}\right)(K_1\oplus\ldots \oplus K_m)\\
\label{12condCtau}
\theta_1K_1\oplus\ldots \oplus \theta_mK_m\subset &C&\subset
\left(1+c_0^n\left(\frac{\varepsilon}{\tau}\right)^{\frac1{95n}}\right)(\theta_1K_1\oplus\ldots \oplus \theta_mK_m).
\end{eqnarray}

Finally, if $\lambda\in[\tau,1-\tau]$ holds for some $\tau\in(0,\frac12]$
in Theorem~\ref{logBMuncondstab} and $\varepsilon\geq \tau$, then choosing again
$m=n$, $K_k=\frac1n(\R e_k\cap K)$, and $\theta_k>0$ in a way such that
$\theta_k(\R e_k\cap K)=\R e_k\cap C$ for $k=1,\ldots,n$, \eqref{approximate-box} yields
\eqref{12condKtau} and \eqref{12condCtau}.
\proofbox

\section{Convex bodies and simplicial cones}
\label{seclinear-reflections}

In this section, we consider the part of a convex body in a Weyl chamber. For a convex body $M$, we write $\partial' M$
to denote the set of every smooth boundary point $x\in \partial M$ where only one unique exterior normal $\nu_{M,x}$ exists, and hence the $(n-1)$-dimensional Hausdorff measure of  $\partial M\backslash \partial' M$ is zero (see Schneider \cite{Sch14}). We recall that the linear $(n-1)$-dimensional subspaces  $H_1,\ldots,H_n\subset\R^n$ are called independent if $H_1\cap\ldots\cap H_n=\{o\}$.

\begin{lemma}
\label{normalsincone}
Let $H_1,\ldots,H_n\subset\R^n$ be independent linear $(n-1)$-dimensional subspaces, and let
$W$ be the closure of a connected component of $\R^n\backslash(H_1\cup\ldots\cup H_n)$.
\begin{description}
\item{(i)}
If $M$ is a convex body in $\R^n$ symmetric through $H_1,\ldots,H_n$, then
$\nu_{M,q}\in W$ for any $q\in W\cap \partial' M$, and in turn
$$
M\cap W=\{x\in W:\,\langle x,u\rangle \leq h_M(u)\;\forall u\in W\}.
$$
\item{(ii)} If $\lambda\in(0,1)$ and $K$ and $C$ are convex bodies in $\R^n$ symmetric through $H_1,\ldots,H_n$, then
$$
W\cap((1-\lambda)K+_0\lambda\,C)=\{x\in W:\,\langle x,u\rangle \leq h_K(u)^{1-\lambda}h_C(u)^\lambda\;\forall u\in W\}.
$$
\end{description}
\end{lemma}
\proof For (i), it is sufficient to prove the first statement; namely, if $q\in {\rm int}\,W\cap \partial' M$, then
$\nu_{M,q}\in W$.

 Let $u_i\in S^{n-1}$, $i=1,\ldots,n$, such that $W=\bigcap_{i=1}^n\{x\in\R^n:\,\langle x,u_i\rangle\geq 0\}$, and hence $u_i$ is a normal to $H_i$ and
$\langle q,u_i\rangle>0$ for $i=1,\ldots,n$, and (i)  is equivalent with the statement that
if $i=1,\ldots,n$, then
\begin{equation}
\label{normalsincone-ui}
\langle u_i,\nu_{M,q}\rangle\geq 0.
\end{equation}
Since $q'=q-2\langle q,u_i\rangle u_i$ is the reflected image of $q$ through $H_i$, we have $q'\in M$; therefore,
$$
0\leq \langle \nu_{M,q},q-q'\rangle=\langle \nu_{M,q},2\langle q,u_i\rangle u_i\rangle=
2\langle q,u_i\rangle\cdot\langle \nu_{M,q}, u_i\rangle.
$$
As $\langle q,u_i\rangle>0$, we conclude \eqref{normalsincone-ui}, and in turn (i).

For (ii), let $M=(1-\lambda)K+_0\lambda\,C$, and let
$$
M_+=\{x\in W:\,\langle x,u\rangle \leq h_K(u)^{1-\lambda}h_C(u)^\lambda\;\forall u\in W\}.
$$
Readily, $W\cap M\subset M_+$. Therefore, (ii) follows if for any $q\in \partial'M\cap{\rm int}W$, we have
$q\in \partial M_+$. As $q\in \partial M\cap{\rm int}W$, there exists $u\in S^{n-1}$ such that
$\langle q,u\rangle= h_K(u)^{1-\lambda}h_C(u)^\lambda$.
Since $q\in \partial' M\cap W$, we have $u=\nu_{M,q}$, and hence (i) yields that $\nu_{M,q}\in W$.
Therefore $q\in\partial M_+$,
proving Lemma~\ref{normalsincone} (ii).
\proofbox

In order to use the known results about unconditional convex bodies, the main idea
is to linearly transfer a Weyl chamber $W$ into the corner $\R_{\geq 0}^n$. For a matrix $\Phi\in{\rm GL}(n,\R)$, its transpose is denoted by $\Phi^\top$, and the inverse of the transpose by $\Phi^{-\top}$.

\begin{lemma}
\label{viuncond}
Let $K$ be a convex body in $\R^n$ with $o\in{\rm int}\,K$, let independent $v_1,\ldots,v_n\in\R^n$
satisfy that $\langle v_i,v_j\rangle\geq 0$ for $1\leq i\leq j\leq n$, let
$W={\rm pos}\,\{v_1,\ldots,v_n\}$, and let $\Phi W=\R_{\geq 0}^n$ for a $\Phi\in {\rm GL}\,(n,\R)$.
\begin{description}
\item{(i)} $\Phi^{-\top} W\subset \R_{\geq 0}^n$.

\item{(ii)} If $\nu_{K,x}\in W$ for all $x\in W\cap \partial' K$, then
\begin{equation}
\label{Phixnormal0}
\nu_{\Phi K,z}\in  \R_{\geq 0}^n\mbox{ \ for all $z\in  \R_{\geq 0}^n\cap \partial' \Phi K$};
\end{equation}

\item{(iii)} and there exists an unconditional convex body $K_0$ such that
$$
\R_{\geq 0}^n\cap K_0=\Phi (W\cap K).
$$
\end{description}
\end{lemma}
\proof Let $e_1,\ldots,e_n$ be the standard orthonormal basis of $\R^n$ indexed in a way such that $e_i=\Phi v_i$.
First we claim that
\begin{equation}
\label{Phixnormal}
\langle e_i,\Phi^{-\top} v\rangle\geq 0 \mbox{ \ for $v\in W$ and $i=1,\ldots,n$.}
\end{equation}
Since $v=\sum_{j=1}^n\lambda_j v_j$ for
$\lambda_1,\ldots,\lambda_n\geq 0$, we deduce from $\langle v_j,v_i\rangle\geq 0$ that
$$
0\leq\left\langle \sum_{j=1}^n\lambda_j v_j,v_i\right\rangle=\langle v,v_i\rangle=
\langle  \Phi^{-\top} v,\Phi v_i\rangle=\langle  \Phi^{-\top} v,e_i\rangle,
$$
proving \eqref{Phixnormal}. In turn, we deduce (i) from \eqref{Phixnormal}.

If $z\in  W\cap \partial' K$, then $\nu_{K, z}\in W$ and $\Phi^{-\top}\nu_{K, z}$ is an exterior normal to $\Phi K$ at $\Phi z$,
therefore, (ii) follows from (i).

Now \eqref{Phixnormal0} yields that if $z=(z_1,\ldots,z_n)\in \R_{\geq 0}^n\cap \partial' \Phi K$
and $0\leq y_i\leq z_i$, $i=1,\ldots,n$, then $y=(y_1,\ldots,y_n)\in \Phi K$. Therefore repeatedly reflecting
$\R_{\geq 0}^n\cap \Phi K$ through the coordinate hyperplanes, we obtain the unconditional convex body $K_0$ such that
$\R_{\geq 0}^n\cap K_0=\R_{\geq 0}^n\cap \Phi K=\Phi (W\cap K)$.
\proofbox

\section{Some properties of Coxeter groups}
\label{secCoxeter-groups}

Since if a  linear map $A$ leaves a convex body $K$ invariant, then the minimal volume Loewner ellipsoid is also invariant under $A$, Barthe, Fradelizi \cite{BaF13} prove that it is sufficient to consider orthogonal reflections in our setting.

\begin{lemma}[Barthe, Fradelizi]
\label{ref-orth}
If the convex bodies $K$ and $C$  in $\R^n$ are invariant under linear reflections $A_1,\ldots,A_n$ through $n$ independent
linear $(n-1)$-planes $H_1,\ldots,H_n$, then there exists $B\in{\rm SL}(n)$ such that
$B A_1 B^{-1},\ldots,B A_n B^{-1}$ are orthogonal reflections through
$B H_1,\ldots,B H_n$ and leave $B K$ and $B C$ invariant.
\end{lemma}

For the theory of Coxeter groups, we follow Humpreys \cite{Hum90}.
For an $n$-dimensional real vector space $V$ equipped with a Euclidean structure, let  $G$ be closure of the Coxeter group generated by the orthogonal reflections through
$p_1^\bot,\ldots,p_n^\bot$ for independent $p_1,\ldots,p_n\in V$.
A linear subspace $L$ of $V$ is invariant under $G$ if and only if $p_1,\ldots,p_n\in L\cup L^\bot$.
We say that an invariant linear subspace $L$ is irreducible if $L\neq\{o\}$ and any invariant subspace
$L'\subset L$ satisfies either $L'=L$ or $L'=\{o\}$, and hence
the action of $G$ on an irreducible invariant subspace is irreducible. Since the intersection and the orthogonal complement of invariant subspaces is invariant, the irreducible subspaces $L_1,\ldots,L_m$, $m\geq 1$ are pairwise orthogonal, and
\begin{equation}
\label{irred-direct-sum}
L_1\oplus \ldots \oplus L_m=V.
\end{equation}
It follows that any $A\in G$ can be written as $A=A|_{L_1}\oplus\ldots\oplus A|_{L_m}$ where $A|_{L_i}$ is the restriction of $A$ to $L_i$ for $i=1,\ldots,m$. For an invariant subspace
$L\subset V$, we set $G|_L=\{A|_L:\, A\in G\}$, and write $O(L)$ to denote the group of isometries of $L$ fixing the origin.
In particular, our main task is to understand irreducible Coxeter groups.

\begin{lemma}[Barthe, Fradelizi]
\label{infty-O(L)}
Let  $G$ be closure of the Coxeter group generated by the orthogonal reflections through
$p_1^\bot,\ldots,p_n^\bot$ for independent $p_1,\ldots,p_n\in \R^n$.
If $L\subset \R^n$ is an irreducible invariant subspace, and $G|_L$ is infinite, then
$G|_L=O(L)$.
\end{lemma}

Next, if $L$ is an irreducible invariant $d$-dimensional linear subspace of $V$ with repect to the closure $G$ of a Coxeter group and $G|_L$ is finite, then a more detailed analysis is needed.
To set up the correponding notation, let $G'=G|_L$ be the finite Coxeter group generated by some orthogonal reflections acting on $L$.  Let $H_1,\ldots,H_k\subset L$ be the linear $(d-1)$-dimensional  subspaces such that the reflections
in $G'$ are the ones through $H_1,\ldots,H_k$, and let $u_1,\ldots,u_{2k}\in L\backslash\{o\}$ be a system of roots for $G'$; namely, there are exactly two roots orthogonal to each $H_i$, and these two roots are opposite.
We note that for algebraic purposes, one usually normalizes the roots in a way such that
$\frac{2\langle u_i,u_j\rangle}{2\langle u_i,u_i\rangle}$ is an integer but we drop this condition because we are only interested in the cones determined by the roots.

Let $W$ be the closure of a Weyl chamber; namely, a connected component of $L\backslash(H_1\cup\ldots\cup H_k)$.
It is known (see \cite{Hum90}) that
$$
W={\rm pos}\{v_1,\ldots,v_d\}=\left\{\sum_{i=1}^d\lambda_iv_i:\,\forall\,\lambda_i\geq 0\right\}
$$
where $v_1,\ldots,v_d\in L$ are independent.
 In addition, for any $x\in L\backslash (H_1\cup\ldots\cup H_k)$, there exists a unique $A\in G'$ such that $x\in AW$, and hence the Weyl chambers are in a natural bijective correspondence with $G'$. We may reindex
$H_1,\ldots,H_k$ and $u_1,\ldots,u_{2k}$ in a way such that $H_i=u_i^\bot$ for $i=1,\ldots,d$ are the "walls" of
$W$, and
\begin{equation}
\label{Weyl-chambers-equation}
\begin{array}{rcll}
\langle u_i,v_i\rangle &>&0&\mbox{ for $i=1,\ldots,d$};\\
\langle u_i,v_j\rangle&=&0&\mbox{ for $1\leq i< j\leq d$.}
\end{array}
\end{equation}
In this case, reflections $L\to L$ through $H_1,\ldots,H_d$ generate $G'$, and
$u_1,\ldots,u_d$ is called a simple system of roots. The order we list simple roots is not related to the corresponding Dynkin diagram.

\begin{lemma}
\label{Weyl-angle}
Let $G$ be the Coxeter group
generated by the orthogonal reflections through
$p_1^\bot,\ldots,p_n^\bot$ for independent $p_1,\ldots,p_n\in \R^n$.
If $L\subset \R^n$ is an irreducible invariant $d$-dimensional subspace with $d\geq 2$, and $G|_L$ is finite,
and $W={\rm pos}\{v_1,\ldots,v_d\}\subset L$ is the closure
of a Weyl chamber  for $G|_L$, then
\begin{equation}
\label{Weyl-angle-eq}
\langle v_i,v_j\rangle\geq  \frac1{d}\cdot \|v_i\|\cdot\|v_j\|.
\end{equation}
\end{lemma}
\proof Let $G'=G|_L$.
We use the classification of finite irreducible Coxeter groups. For the cases when $G'$ is either of
$D_d,E_6,E_7,E_8$ (see Adams \cite{Ada16} about $E_6,E_7,E_8$), we use the known simple systems of roots in terms of the orthonormal basis
$e_1,\ldots,e_d$ of $L$ to construct $v_1,\ldots,v_d$ {\it via} \eqref{Weyl-chambers-equation}.
However, there is a unified construction for the other finite irreducible Coxeter groups because they are the symmetries of some regular polytopes.   \\

\noindent{\bf Case 1: $G'$ is one of the types} $\mathbf{I_2(m)}$, $\mathbf{A_d}$, $\mathbf{B_d}$,
$\mathbf{F_4}$, $\mathbf{H_3}$, $\mathbf{H_4}$\\
In this case, $G'$ is the symmetry group of some $d$-dimensional regular polytope $P$ centered at the origin.
Let $F_0\subset \ldots\subset F_{d-1}$ be a tower of faces of $P$ where ${\rm dim}\,F_i=i$, $i=0,\ldots,d-1$. Defining
$v_i$ to be the centroid of $F_{i-1}$, $i=1,\ldots,d$, we have that
$W={\rm pos}\{v_1,\ldots,v_d\}$ is the closure of a Weyl chamber because the symmetry group of $P$ is simply transitive on the towers of faces of $P$.

As $G'$ is irreducible, the John ellipsoid of $P$ (the unique ellipsoid of largest volume contained in $P$) is a $d$-dimensional
ball centered at the origin of some radius $r>0$. It follows that $P\subset dr B^d$, and hence
$r\leq \|v_i\|\leq dr$ for $i=1,\ldots,d$. In addition, $v_i$ is the closest point of ${\rm aff}\,F_{i-1}$ to the origin
for $i=1,\ldots,d$, and $v_j\in F_{i-1}$ if $1\leq j\leq i$,  thus $\langle v_j,v_i\rangle=\langle v_i,v_i\rangle$
if $1\leq j\leq i\leq d$. We conclude that if $1\leq j\leq i\leq d$, then
$$
\frac{\langle v_j,v_i\rangle}{\|v_j\|\cdot\|v_i\|}=
\frac{\|v_i\|}{\|v_j\|}\geq \frac1d.
$$
\mbox{ }\\
\noindent{\bf Case 2: $G'=D_n$}\\
In this case, a simple system of roots is
$$
\begin{array}{rcll}
u_i&=&e_i-e_{i+1}&\mbox{ for $i=1,\ldots,d-1$},\\
u_d&=&e_{d-1}+e_d.&
\end{array}
$$
In turn, we may choose $v_1,\ldots,v_d$ as
$$
\begin{array}{rcll}
v_i&=&\sum_{l=1}^ie_l&\mbox{ for $i=1,\ldots,d-2$ and $i=d$},\\
v_{d-1}&=&-v_d+\sum_{l=1}^{d-1}e_l.&
\end{array}
$$
As $\langle v_i,v_j\rangle$ is a positive integer for $i\neq j$,
and $\|v_i\|\leq \sqrt{d}$ for $i=1,\ldots,d$, we conclude (\ref{Weyl-angle-eq}).\\

\noindent{\bf Case 3: $G'=E_6$}\\
In this case $d=6$, and a simple system of roots is
$$
\begin{array}{rcll}
u_i&=&e_i-e_{i+1}&\mbox{ for $i=1,2,3,4$},\\
u_5&=&e_4+e_5&\\
u_6&=&\sqrt{3}\,e_6-\sum_{l=1}^5e_l.&
\end{array}
$$
Using coordinates in $e_1,\ldots,e_6$, we may choose $v_1,\ldots,v_6$ as
$v_1=(\sqrt{3},0,0,0,0,1)$,
$v_2=(\sqrt{3},\sqrt{3},0,0,0,2)$,
$v_3=(\sqrt{3},\sqrt{3},\sqrt{3},0,0,3)$,
$v_4=(1,1,1,1,-1,\sqrt{3})$,
$v_5=(1,1,1,1,1,\frac{5}{\sqrt{3}})$ and
$v_6=(0,0,0,0,0,3)$.
As $\langle v_i,v_j\rangle\geq 3$ for $i\neq j$,
and $\|v_i\|\leq \sqrt{18}$ for $i=1,\ldots,6$, we conclude (\ref{Weyl-angle-eq}).\\

\noindent{\bf Case 4: $G'=E_7$}\\
In this case $d=7$, and a simple system of roots is
$$
\begin{array}{rcll}
u_i&=&e_i-e_{i+1}&\mbox{ for $i=1,2,3,4,5$},\\
u_6&=&e_5+e_6&\\
u_7&=&\sqrt{2}\,e_7-\sum_{l=1}^6e_l.&
\end{array}
$$
Using coordinates in $e_1,\ldots,e_7$, we may choose $v_1,\ldots,v_7$ as
$v_1=(2,0,0,0,0,0,\sqrt{2})$,
$v_2=(1,1,0,0,0,0,\sqrt{2})$,
$v_3=(1,1,1,0,0,0,\frac{3}{\sqrt{2}})$,
$v_4=(1,1,1,1,0,0,2\sqrt{2})$,
$v_5=(1,1,1,1,1,-1,2\sqrt{2})$,
$v_6=(1,1,1,1,1,1,3\sqrt{2})$ and
$v_7=(0,0,0,0,0,4)$.
As $\langle v_i,v_j\rangle\geq 4$ for $i\neq j$,
and $\|v_i\|< \sqrt{28}$ for $i=1,\ldots,7$, we conclude (\ref{Weyl-angle-eq}).\\

\noindent{\bf Case 5: $G'=E_8$}\\
In this case $d=8$, and a simple system of roots is
$$
\begin{array}{rcll}
u_i&=&e_i-e_{i+1}&\mbox{ for $i=1,2,3,4,5,6,7$},\\
u_8&=&-\sum_{l=1}^5e_l+\sum_{l=6}^8e_l.&
\end{array}
$$
Using coordinates in $e_1,\ldots,e_8$, we may choose $v_1,\ldots,v_8$ as
$v_1=(1,-1,-1,-1,-1,-1,-1,-1)$,
$v_2=(0,0,-1,-1,-1,-1,-1,-1)$,
$v_3=(-1,-1,-1,-3,-3,-3,-3,-3)$,
$v_4=(-1,-1,-1,-1,-2,-2,-2,-2)$,
$v_5=(-1,-1,-1,-1,-1,-\frac{5}3,-\frac{5}3,-\frac{5}3)$,
$v_6=(-1,-1,-1,-1,-1,-1,-2,-2)$,
$v_7=(-1,-1,-1,-1,-1,-1,-1,-3)$ and
$v_8=(-1,-1,-1,-1,-1,-1,-1,-1)$.
As $\langle v_i,v_j\rangle\geq 6$ for $i\neq j$,
and $\|v_i\|< \sqrt{48}$ for $i=1,\ldots,8$, we conclude (\ref{Weyl-angle-eq}).
\proofbox

For a convex body invariant under a Coxeter group, we can determine some exterior normal at certain points provided by the symmetries of the convex body.

\begin{lemma}
\label{normals-symmetry}
Let $G$ be the closure of a Coxeter group generated by $n$ independent orthogonal reflections of $\R^n$,
let $L\subset\R^n$ be an irreducible linear subspace
and let $K$ be a convex body in $\R^n$ invariant under $G$.
\begin{description}
\item{(i)} If $G|_L$ is finite, and $W={\rm pos}\{v_1,\ldots,v_d\}\subset L$ is the closure
of a Weyl chamber  for $G|_L$, and $t_iv_i\in\partial K$ for $t_i>0$, $i=1,\ldots,d$, then
$v_i$ is an exterior normal at $tv_i$.
\item{(ii)} If $G|_L$ is infinite and $v\in L\backslash \{o\}$, and $tv\in\partial K$ for $t>0$, then
$v$ is an exterior normal at $tv$.
\end{description}
\end{lemma}
\proof Let $d={\rm dim}\,L$.

For (i), first we claim that there exist independent $u_1,\ldots,u_{n-1}\in v_i^\bot$ such that the reflection through
$u_j^\bot$ lies in $G$ for $j=1,\ldots,n-1$.
To construct $u_1,\ldots,u_{n-1}\in v_i^\bot$, if $d\geq 2$, then we choose roots $u_1,\ldots,u_{d-1}\in v_i^\bot$ for
$G|_L$ that corresponds to the walls of $W$ containing $v_i$. In addition, if $d<n$, then
 we choose independent  $u_d,\ldots,u_{n-1}\in L^\bot$ such that the reflection through
$u_j^\bot$ lies in $G$ for $j=d,\ldots,n-1$, completing the construction of $u_1,\ldots,u_{n-1}$.

Let $N=\{y\in\R^n:\,\langle y,t_iv_i-x\rangle\geq 0\;\forall x\in K\}$ be the normal cone at $t_iv_i\in\partial K$.
If $N=\R_{\geq 0}v_i$, then we are done; therefore, we assume that $N\neq\R_{\geq 0}v_i$. Since $o\in{\rm int}\,K$, 
$\langle y,v_i\rangle>0$ for any $y\in N\backslash\{o\}$, and since 
$N$ is a cone and $N\neq\R_{\geq 0}v_i$,
there exists $w\in v_i^\bot\backslash\{o\}$ such that $z=v_i+w\in N$. Let $H\subset G$ be the closure of the subgroup generated by the reflections through $u_1^\bot,\ldots,u_{n-1}^\bot$, and hence both $\R v_i$ and $v_i^\bot$ are invariant under  $H$.
Since $u_1,\ldots,u_{n-1}\in v_i^\bot$ are independent, the centroid of
$M={\rm conv}\{Aw:\,A\in H\}\subset v_i^\bot$ is $o$. We deduce that
the centroid of $v_i+M={\rm conv}\{Az:\,A\in H\}\subset N$ is $v_i$; therefore, $v_i\in N$.

For (ii), the argument is essentially same because similarly,  there exist independent
$\tilde{u}_1,\ldots,\tilde{u}_{n-1}\in v^\bot$ such that the reflection through
$\tilde{u}_j^\bot$ lies in $G$ for $j=1,\ldots,n-1$.
\proofbox

\section{The proof Theorem~\ref{logBMstab} }
\label{seclogBMstab}

Lemma~\ref{ref-orth} and the linear invariance of the $L_0$-sum yield that we may assume that $A_1,\ldots,A_n$ are orthogonal reflections through the linear $(n-1)$-spaces $H_1,\ldots,H_n$, respectively,  with $H_1\cap\ldots\cap H_n=\{o\}$
where $K$ and $C$ are invariant under $A_1,\ldots,A_n$.

Let $G$ be the closure of the group generated by
$A_1,\ldots,A_n$, and let $L_1,\ldots,L_m$ be the irreducible invariant subspaces of $\R^n$ of the action of $G$.
If $t_1,\ldots,t_m>0$ and $\Psi\in{\rm GL}(n,\R)$   satisfies $\Psi x=t_ix$ for $x\in L_i$ and $i=1,\ldots, m$, then
\begin{equation}
\label{diag-irredtiLi}
\mbox{$\Psi K$ and $\Psi C$ are both invariant under $G$.}
\end{equation}

Let $E$ be the John ellipsoid of $K$, that is, the unique ellipsoid of maximal volume contained in $K$. Therefore, $E$ is also invariant under $G$. In particular, we can choose the principal directions of $E$ in a way such that each is contained in one of the $L_i$, and $L_i\cap E$  is a Euclidean ball of dimension ${\rm dim}\,L_i$. Therefore, after applying a suitable linear transformation like in \eqref{diag-irredtiLi}, we may assume that $E=B^n$, and hence
\begin{equation}
\label{John-ball-inK}
B^n\subset K\subset nB^n.
\end{equation}

For any $i=1,\ldots,n$, let $G_i=G|_{L_i}$ if $G|_{L_i}$ is finite, and let $G_i$ be the symmetry group
of some ${\rm dim}\,L_i$ dimensional regular simplex in $L_i$ centered at the origin if $G|_{L_i}$ is infinite.

We consider the finite subgroup $\widetilde{G}\subset G$ that is the direct sum of $G_1,\ldots,G_m$, acting in the natural way
$\widetilde{G}|_{L_i}=G_i$ for $i=1,\ldots,m$.
Let $0=p_0<p_1<\ldots<p_m=n$ satisfy that
$p_i-p_{i-1}={\rm dim}\,L_i$ for $i=1,\ldots,m$.
We choose a basis $v_1,\ldots,v_n\in S^{n-1}$ of $\R^n$,
 in a way such that
for each $i=1,\ldots,m$, $W_i={\rm pos}\{v_{p_{i-1}+1},\ldots,v_{p_i}\}$ is the closure
of a Weyl chamber for the irreducible action of $G_i$ on $L_i$.

According to Lemma~\ref{Weyl-angle}, these $v_1,\ldots,v_n\in S^{n-1}$ satisfy that
\begin{eqnarray}
\label{vi-angles-big}
\langle v_j,v_l\rangle&\geq& \frac1n\mbox{ \ \ \ if $p_{i-1}+1\leq j<l\leq p_i$ and $i=1,\ldots,m$;}\\
\label{vi-angles-ort}
\langle v_j,v_l\rangle&=& 0\mbox{ \ \ \ if there exists $i=1,\ldots,m-1$ such that $j\leq p_i< l$.}
\end{eqnarray}
 Let $e_1,\ldots,e_n$ be the standard orthonormal basis of $\R^n$,
 let $\Phi\in {\rm GL}(n)$ satisfy that $\Phi v_i=e_i$, $i=1,\ldots,n$, and let
$$
W=W_1\oplus\ldots \oplus W_m.
$$
It follows that $\Phi W=\R_{\geq 0}^n$ and ${\rm int}\,W$ is a fundamental domain for $\widetilde{G}$ in the sense that
\begin{equation}
\label{W-Fund-domain}
\begin{array}{rcll}
\bigcup\{AW:\,A\in\widetilde{G}\}&=&\R^n&\\
{\rm int}\,AW\cap {\rm int}\,BW&=&\emptyset&\mbox{ \ if $A,B\in\widetilde{G}$ and $A\neq B$.}
\end{array}
\end{equation}

If $i\in\{1,\ldots,m\}$ and $p_{i-1}+1\leq j\leq p_i$, then we define $u_j\in L_i\cap S^{n-1}$
by $\langle u_j,v_j\rangle>0$ and $\langle u_j,v_l\rangle=0$ for $l\neq j$. Therefore,
$u_1^\bot,\ldots,u_n^\bot$ are the walls of $W$; namely, the linear hulls of the facets of the simplicial cone $W$,
and the reflections through $u_1^\bot,\ldots,u_n^\bot$ are symmetries of both $K$ and $C$
(and actually generate $\widetilde{G}$).
We may apply Lemma~\ref{viuncond} to $W$ because of Lemma~\ref{normalsincone},
\eqref{vi-angles-big} and \eqref{vi-angles-ort}, and deduce
 the existence of unconditional convex bodies $\widetilde{K}$ and $\widetilde{C}$ such that
$$
\R_{\geq 0}^n\cap \widetilde{K}=\Phi (W\cap K)\mbox{ and }\R_{\geq 0}^n\cap \widetilde{C}=\Phi (W\cap C).
$$

We claim that
\begin{equation}
\label{PhiWL0sum}
\R_{\geq 0}^n\cap((1-\lambda)\widetilde{K}+\lambda \widetilde{C})\subset
\Phi\left(W\cap((1-\lambda)K+_0\lambda C)\right).
\end{equation}
According to Lemma~\ref{normalsincone}  and to
 $\Phi^{-\top} W\subset \R^n_{\geq 0}$ ({\it cf.} Lemma~\ref{viuncond}), we have
\begin{eqnarray*}
\R_{\geq 0}^n\cap((1-\lambda)\widetilde{K}+\lambda \widetilde{C})&=&
\{x\in \R_{\geq 0}^n:\,\langle x,u\rangle\leq h_{\widetilde{K}}(u)^{1-\lambda}h_{\widetilde{C}}(u)^{1-\lambda}\;
\forall u\in\R_{\geq 0}^n\}\\
&\subset&
\{x\in \R_{\geq 0}^n:\,\langle x,u\rangle\leq h_{\widetilde{K}}(u)^{1-\lambda}h_{\widetilde{C}}(u)^\lambda\;\forall u\in\Phi^{-\top} W\}.
\end{eqnarray*}
We observe that if $u\in \Phi^{-\top} W$, then  there exist
$y_0\in\R_{\geq 0}^n\cap\partial \widetilde{K}=\R_{\geq 0}^n\cap\partial (\Phi K)$ and
$z_0\in\R_{\geq 0}^n\cap\partial \widetilde{C}=\R_{\geq 0}^n\cap\partial (\Phi C)$ with
$h_{\widetilde{K}}(u)=\langle y_0,u\rangle$
and $h_{\widetilde{C}}(u)=\langle z_0,u\rangle$. For $v=\Phi^{\top} u\in W$, $y=\Phi^{-1}y_0\in W\cap \partial K$
and $y=\Phi^{-1}y_0\in W\cap \partial K$, it follows that $v$ is an exterior normal to $K$ at $y$ and to $C$ at $z$,
and
$$
h_{\widetilde{K}}(u)^{1-\lambda}h_{\widetilde{C}}(u)^\lambda
=\langle \Phi y,\Phi^{-\top}v\rangle^{1-\lambda}\langle \Phi z,\Phi^{-\top}v\rangle^\lambda
=\langle y,v\rangle^{1-\lambda}\langle z,v\rangle^\lambda=
h_K(v)^{1-\lambda}h_C(v)^\lambda.
$$
We deduce from the considerations just above and from applying Lemma~\ref{normalsincone} to $W$ that
\begin{eqnarray*}
\R_{\geq 0}^n\cap((1-\lambda)\widetilde{K}+\lambda \widetilde{C})
&\subset&
\Phi\{q\in W:\,\langle q,v\rangle\leq h_{K}(v)^{1-\lambda}h_{K}(v)^\lambda\;\forall v\in W\}\\
&=& \Phi\left(W\cap((1-\lambda)K+_0\lambda C)\right),
\end{eqnarray*}
proving \eqref{PhiWL0sum}.

Writing $|\widetilde{G}|$ to denote the cardinality of $\widetilde{G}$, \eqref{W-Fund-domain} yields
$$
V(M)=|\widetilde{G}|\cdot V(M\cap W)
$$
where $M$ is either $K$, $C$ or $(1-\lambda)\cdot K +_0 \lambda\cdot C$.
We deduce from \eqref{PhiWL0sum} and the condition in Theorem~\ref{logBMstab} that
\begin{eqnarray*}
V((1-\lambda)\cdot \widetilde{K} +_0 \lambda\cdot \widetilde{C})&= &
2^nV\left(\R_{\geq 0}^n\cap((1-\lambda)\cdot \widetilde{K} +_0 \lambda\cdot \widetilde{C})\right)\\
&\leq & 2^nV\left(\Phi\left(W\cap((1-\lambda)K+_0\lambda C)\right)\right)\\
&\leq&\frac{2^n|\det\Phi|}{|\widetilde{G}|}\cdot(1+\varepsilon) V(K)^{1-\lambda} V(C)^\lambda\\
&=&(1+\varepsilon) V(\widetilde{K})^{1-\lambda} V(\widetilde{C})^\lambda.
\end{eqnarray*}
We apply the following equivalent form of Theorem~\ref{logBMuncondstab} to
$\widetilde{K}$ and $\widetilde{C}$
where $\lambda\in [\tau,1-\tau]$ for $\tau\in(0,\frac12]$. There exist absolute constant $\tilde{c}>1$,
complementary coordinate linear subspaces $\widetilde{\Lambda}_1,\ldots,\widetilde{\Lambda}_k$, $k\geq 1$,
with $\bigoplus_{j=1}^k\widetilde{\Lambda}_j=\R^n$ such that
\begin{equation}
\label{tildeKjKLambda}
\bigoplus_{j=1}^k\left(\widetilde{K}\cap \widetilde{\Lambda}_j\right)\subset
\left(1+\tilde{c}^n\left(\frac{\varepsilon}{\tau}\right)^{\frac1{95n}}\right) \widetilde{K},
\end{equation}
and there exist $\theta_1,\ldots,\theta_k>0$ such that
\begin{equation}
\label{tildeCjCLambda}
\bigoplus_{j=1}^k\theta_j\left(\widetilde{K}\cap \widetilde{\Lambda}_j\right)\subset
\widetilde{C}\subset
\left(1+\tilde{c}^n\left(\frac{\varepsilon}{\tau}\right)^{\frac1{95n}}\right)
\bigoplus_{j=1}^k\theta_j\left(\widetilde{K}\cap \widetilde{\Lambda}_j\right).
\end{equation}
For $\Lambda_j=\Phi^{-1}\widetilde{\Lambda}_j$, $j=1,\ldots,k$, we deduce that
\begin{equation}
\label{WKjKLambda}
W\cap \sum_{j=1}^k\left(K\cap \Lambda_j\right)\subset
\left(1+\tilde{c}^n\left(\frac{\varepsilon}{\tau}\right)^{\frac1{95n}}\right) (W\cap K),
\end{equation}
and
\begin{equation}
\label{WCjCLambda}
W\cap \sum_{j=1}^k\theta_j\left(K\cap \Lambda_j\right)\subset
W\cap C\subset
\left(1+\tilde{c}^n\left(\frac{\varepsilon}{\tau}\right)^{\frac1{95n}}\right)
\left(W\cap \sum_{j=1}^k\theta_j\left(K\cap \Lambda_j\right)\right).
\end{equation}
We observe that each $\Lambda_j$ is spanned by a subset of $v_1,\ldots,v_n$.

For the rest of the argument, first we assume that $\varepsilon$ is small enough to satisfy
\begin{equation}
\label{Lambdaj-invariant-eps}
\tilde{c}^n\left(\frac{\varepsilon}{\tau}\right)^{\frac1{95n}}<\frac1{n^2}.
\end{equation}

We claim that if \eqref{Lambdaj-invariant-eps} holds, then
\begin{equation}
\label{Lambdaj-invariant}
\mbox{each $\Lambda_j$, $j=1,\ldots,k$, is invariant under $G$.}
\end{equation}
We suppose indirectly that the claim \eqref{Lambdaj-invariant} does not hold, and we seek a contradiction.
In this case, $k\geq 2$. Since each $\Lambda_j$ is spanned by a subset of $v_1,\ldots,v_n$,
after possibly reindexing $L_1,\ldots,L_m$, $\Lambda_1,\ldots,\Lambda_k$ and
 $v_1,\ldots,v_n$, we may assume that
$v_1\in L_1\cap\Lambda_1$ and $v_2\in L_1\cap\Lambda_2$. For $i=1,\ldots,n$, let $s_i>0$ satisfy
$s_iv_i\in\partial K$; therefore, \eqref{John-ball-inK} yields
\begin{equation}
\label{si1n-bound}
1\leq s_i\leq n,
\end{equation}
and hence
\begin{equation}
\label{s12L1Lambda12}
\mbox{$s_1v_1\in L_1\cap K\cap\Lambda_1$ and $v_2\in L_1\cap K\cap\Lambda_2$}.
\end{equation}
It follows from \eqref{vi-angles-big} that
\begin{equation}
\label{v1v2-angles-big}
\langle v_1,v_2\rangle\geq \frac1n.
\end{equation}
We deduce from \eqref{s12L1Lambda12}, and then from \eqref{WKjKLambda} that
\begin{equation}
\label{s1v1+v2}
s_1v_1+v_2\in
W\cap \sum_{j=1}^k\left(K\cap \Lambda_j\right)\subset
\left(1+\tilde{c}^n\left(\frac{\varepsilon}{\tau}\right)^{\frac1{95n}}\right) (W\cap K).
\end{equation}
Lemma~\ref{normals-symmetry} yields that $v_1$ is an exterior unit normal to $\partial K$ at $s_1v_1$,
and hence $s_1=h_K(v_1)$. We deduce from first \eqref{s1v1+v2} and then from
assumption \eqref{Lambdaj-invariant-eps} and the formula \eqref{si1n-bound} that
\begin{eqnarray}
\nonumber
s_1+\langle v_1,v_2\rangle&=&\langle v_1,s_1v_1+v_2\rangle
\leq\left(1+\tilde{c}^n\left(\frac{\varepsilon}{\tau}\right)^{\frac1{95n}}\right)h_K(v_1)\\
\label{s1v1+v2proj-up}
&=&s_1+\tilde{c}^n\left(\frac{\varepsilon}{\tau}\right)^{\frac1{95n}}s_1<s_1+\frac1n.
\end{eqnarray}
On the other hand, we have $s_1+\langle v_1,v_2\rangle\geq s_1+\frac1n$ by \eqref{v1v2-angles-big},
contradicting \eqref{s1v1+v2proj-up}. In turn, we conclude \eqref{Lambdaj-invariant}
under the assumption \eqref{Lambdaj-invariant-eps}.

We deduce from \eqref{WKjKLambda}, \eqref{WCjCLambda}, \eqref{Lambdaj-invariant} and the symmetries of $K$ and $C$ that
\begin{equation}
\label{KjKLambda}
\bigoplus_{j=1}^k\left(K\cap \Lambda_j\right)\subset
\left(1+\tilde{c}^n\left(\frac{\varepsilon}{\tau}\right)^{\frac1{95n}}\right) K,
\end{equation}
and
\begin{equation}
\label{CjCLambda}
\bigoplus_{j=1}^k\theta_j\left(K\cap \Lambda_j\right)\subset
C\subset
\left(1+\tilde{c}^n\left(\frac{\varepsilon}{\tau}\right)^{\frac1{95n}}\right)
 \bigoplus_{j=1}^k\theta_j\left(K\cap \Lambda_j\right).
\end{equation}
In addition, the symmetries of $K$ and \eqref{Lambdaj-invariant} yield that $K\cap \Lambda_j=K|\Lambda_j$ for
$j=1,\ldots,k$, therefore,
$$
K\subset \bigoplus_{j=1}^k\left(K\cap \Lambda_j\right).
$$
Combining this relation with \eqref{KjKLambda} and \eqref{CjCLambda} implies
Theorem~\ref{logBMstab}
under the assumption \eqref{Lambdaj-invariant-eps}.

Finally, we assume that
\begin{equation}
\label{Lambdaj-invariant-epsup}
\tilde{c}^n\left(\frac{\varepsilon}{\tau}\right)^{\frac1{95n}}\geq\frac1{n^2},
\end{equation}
and hence
\begin{equation}
\label{Lambdaj-invariant-epsup0}
(5\tilde{c})^n\left(\frac{\varepsilon}{\tau}\right)^{\frac1{95n}}\geq n^2.
\end{equation}

For $i=1,\ldots,m$, the symmetries of $K$ and $C$ yield that $r_i(B^n\cap L_i)$ is the John ellipsoid
of $K\cap L_i$ and $\theta_ir_i(B^n\cap L_i)$ is the John ellipsoid
of $C\cap L_i$ for some $r_i,\theta_i>0$.
For $K_i=\frac{r_i}n\,(B^n\cap L_i)$, $i=1,\ldots,m$, we have
$$
\bigoplus_{i=1}^m K_i\subset{\rm conv}\{mK_1,\ldots,mK_m\};
$$
therefore,
it follows from \eqref{Lambdaj-invariant-epsup0}
that
\begin{eqnarray*}
\bigoplus_{i=1}^m K_i &\subset K\subset&n^2\cdot \bigoplus_{i=1}^m K_i\subset
\left(1+(5\tilde{c})^n\left(\frac{\varepsilon}{\tau}\right)^{\frac1{95n}}\right)\bigoplus_{i=1}^m K_i\\
\bigoplus_{i=1}^m \theta_iK_i&\subset C\subset&n^2\cdot \bigoplus_{i=1}^m \theta_iK_i\subset
\left(1+(5\tilde{c})^n\left(\frac{\varepsilon}{\tau}\right)^{\frac1{95n}}\right)\bigoplus_{i=1}^m \theta_iK_i,
\end{eqnarray*}
proving Theorem~\ref{logBMstab}
under the assumption \eqref{Lambdaj-invariant-epsup}.
\proofbox

\section{Proof of Theorem~\ref{logMstab}}
\label{seclogMstab}

As in the case of Theorem~\ref{logBMstab}, it follows from Lemma~\ref{ref-orth} and the linear invariance of the $L_0$-sum that we may assume that
$A_1,\ldots,A_n$ are orthogonal reflections through the linear $(n-1)$-spaces $H_1,\ldots,H_n$, respectively,  with $H_1\cap\ldots\cap H_n=\{o\}$
where $K$ and $C$ are invariant under $A_1,\ldots,A_n$.
We write $G$ to denote the closure of the group generated by
$A_1,\ldots,A_n$, and $L_1,\ldots,L_m$ to denote the irreducible invariant subspaces of $\R^n$ of the action of $G$.

 For the logarithmic Minkowski Conjecture~\ref{logMconj},
replacing either $K$ or $C$ by a dilate does not change the difference of the two sides; therefore, we may assume that
$$
V(K)=V(C)=1.
$$
In this case, the condition in Theorem~\ref{logMstab} states that
\begin{equation}
\label{logMstab-cond}
\int_{S^{n-1}}\log \frac{h_C}{h_K}\,dV_K<\varepsilon
\end{equation}
for $\varepsilon>0$.

First we assume that
\begin{equation}
\label{neps}
n\varepsilon<1,
\end{equation}
for $t\in[0,1]$, we define
$$
\varphi(t)=V((1-t)\cdot K +_0 t\cdot C).
$$
According to (3.7) in B\"or\"oczky, Lutwak, Yang, Zhang \cite{BLYZ12}, we have
\begin{equation}
\label{logminkder}
\varphi'(0)=n\int_{S^{n-1}}\log \frac{h_C}{h_K}\,dV_K,
\end{equation}
and hence \eqref{logMstab-cond} and the assumption \eqref{neps} yield that
$\varphi'(0)<n\varepsilon$ where $n\varepsilon<1$.
We deduce from Lemma~\ref{log-concave-error} (ii) that
$$
V\left(\frac12\cdot K +_0 \frac12\cdot C\right)=\varphi\left(\frac12\right)<1+n\varepsilon.
$$
Now we apply Theorem~\ref{logBMstab}, and conclude that
 for some $m\geq 1$,
there exist $\theta_1,\ldots,\theta_m>0$ and compact convex sets
$K_1,\ldots,K_m>0$ invariant under $G$ such that ${\rm lin}\,K_i$, $i=1,\ldots,m$, are complementary coordinate subspaces, and
\begin{eqnarray}
\label{logMstab-smalleps1}
K_1\oplus\ldots \oplus K_m\subset &K&\subset
\left(1+c^n\varepsilon^{\frac1{95n}}\right)(K_1\oplus\ldots \oplus K_m)\\
\label{logMstab-smalleps2}
\theta_1K_1\oplus\ldots \oplus \theta_mK_m\subset &C&\subset
\left(1+c^n\varepsilon^{\frac1{95n}}\right)(\theta_1K_1\oplus\ldots \oplus \theta_mK_m)
\end{eqnarray}
where $c>1$ is an absolute constant. In turn, we deduce Theorem~\ref{logMstab} under the assumption $n\varepsilon<1$ on
\eqref{neps}.

On the other hand, if $n\varepsilon\geq 1$, then Theorem~\ref{logMstab} can be proved as
Theorem~\ref{logBMstab}
under the assumption \eqref{Lambdaj-invariant-epsup}.
\proofbox

\noindent{\bf Acknowledgement } We would like to thank Gaoyong Zhang and Richard Gardner for illuminating discussions, and the referees who spotted many typos and mistakes in an earlier version of the paper, and whose remarks have substantially improved the readability of the paper.

\end{document}